\documentclass[11pt]{amsart}

\newcommand{\C}{{\mathbf C}}
\newcommand{\F}{{\mathbf F}}
\newcommand{\Z}{{\mathbf Z}}

\newcommand{\PP}{{\mathbf P}}
\newcommand{\CC}{{\mathcal C}}
\newcommand{\CF}{{\mathcal F}}
\newcommand{\Eb}{{\overline{E}}}
\newcommand{\Fb}{{\overline{F}}}
\newcommand{\fb}{{\overline{f}}}
\newcommand{\gb}{{\overline{g}}}
\newcommand{\Q}{{\mathbf Q}}
\newcommand{\Qbar}{{\overline{\Q}}}

\newcommand{\kbar}{\overline{k}}
\newcommand{\can}{{\text{can}}}
\newcommand{\charact}{\operatorname{char}}

\newcommand{\Aut}{\operatorname{Aut}}
\newcommand{\Gal}{\operatorname{Gal}}
\newcommand{\Hom}{\operatorname{Hom}}
\newcommand{\Norm}{\operatorname{Norm}}
\newcommand{\Jac}{\operatorname{Jac}}
\newcommand{\SL}{\operatorname{SL}}
\newcommand{\Spec}{\operatorname{Spec}}
\newcommand{\End}{\operatorname{End}}

\newcommand{\isom}{\cong}

\newcommand\ra{\rightarrow}
\newcommand\mapright[1]{\mathop{\longrightarrow}\limits^{#1}}
\newcommand\Bmapdown[1]{\vbox{\vbox to 4pt{}\vbox{\hbox{
             \Big\downarrow\rlap{$\vcenter{\hbox{$\scriptstyle#1$}}$}}\vfill}}}
\newcommand\Bmapup[1]  {\vbox{\vbox to 4pt{}\vbox{\hbox{
               \Big\uparrow\rlap{$\vcenter{\hbox{$\scriptstyle#1$}}$}}\vfill}}}

\newcommand\buff{\rule[-.1ex]{0pt}{2.3ex}}

\newcommand\bbuff{\rule[-2.4ex]{0pt}{6.1ex}}

\newtheorem{theorem}{Theorem}
\newtheorem{lemma}[theorem]{Lemma}
\newtheorem{cor}[theorem]{Corollary}
\newtheorem{prop}[theorem]{Proposition}

\theoremstyle{definition}

\theoremstyle{remark}
\newtheorem{rem}{Remark}	
	
\newtheorem{example}{Example}	

\begin{document}

\title[Torsion subgroups of Jacobians]{Large torsion subgroups of split Jacobians of curves of genus two or three}
\subjclass{Primary 14G05; Secondary 11G30, 14H25, 14H45}
\keywords{Torsion, Jacobian, elliptic differential.}

\author{Everett W.\ Howe}
\address{Center for Communications Research, 4320 Westerra Court, San Diego, CA 92121-1967, USA.}
\email{however@alumni.caltech.edu}

\author{Franck Lepr\'{e}vost}
\address{Universit\'e Pierre et Marie Curie,
CNRS- UMR 9994 Arithmetique,
Tour 46-56, 5\`eme \'etage Case 247; 
4, place Jussieu; 
F-75252 Paris Cedex 05; France.
and 
Technische Universit\"at Berlin,
Fachbereich Mathematik MA 8-1; 
Stra{\ss}e des 17. Juni 136;
D-10623 Berlin; Germany.}
\email{leprevot@math.tu-berlin.de,  leprevot@math.jussieu.fr}

\author{Bjorn Poonen}
\thanks{The first author is supported partially
by NSA Young Investigator Grant MDA904-95-H-1044.
The second author thanks the Max-Planck-Institut f\"ur Mathematik
of Bonn for its hospitality.
The third author is supported by an NSF Mathematical Sciences
Postdoctoral Research Fellowship.
The final version of this paper has appeared in
{\em Forum Math.}\ {\bf 12} (2000), 315--364.
}
\address{Department of Mathematics,
University of California,
Berkeley, CA 94720-3840, USA.}
\email{poonen@math.berkeley.edu}
\date{13 February 1997; revised 19 September 1998.}

\begin{abstract}
We construct examples of families of curves of genus~2 or~3 over $\Q$ whose
Jacobians split completely and have various large rational torsion subgroups.
For example, the rational points on a certain elliptic surface
over $\PP^1$ of positive rank parameterize a family of genus-2 curves over $\Q$
whose Jacobians each have $128$ rational torsion points.
Also, we find the genus-3 curve
$$15625(X^4 + Y^4 + Z^4) - 96914(X^2 Y^2 + X^2 Z^2 + Y^2 Z^2) = 0, \text{\phantom{mmmmmmmmm}}$$
whose Jacobian has $864$ rational torsion points.
\end{abstract}

\maketitle

\section{Introduction}
\label{introduction}

Nearly twenty years ago Mazur settled the question of which groups can occur
as the group of rational torsion points on an elliptic curve over $\Q$,
but the analogous question for Jacobian varieties of curves over $\Q$ of genus 
greater than $1$ remains open.
Most of the work that has been done on this question has centered on the problem
of finding groups that {\em do\/} occur as rational torsion subgroups of Jacobians.
Several researchers have produced families of genus-$2$ curves whose Jacobians
contain various given groups in their rational torsion 
(see \cite{leprevost1}, \cite{leprevost2}, \cite{leprevost5}, \cite{leprevost6},
\cite{ogawa}, and the summary in \cite{leprevost9})
while others have constructed families of curves in which the size of the rational
torsion subgroup of the Jacobian increases as the genus of the curve increases
(see \cite{flynnlarge}, \cite{flynnsequences}, \cite{leprevost3}, \cite{leprevost4},
\cite{leprevost7}, \cite{leprevost8}).
The largest group of rational torsion heretofore known to exist on the Jacobian
of a curve of genus~$2$ was a group of order $30$; for genus-$3$ curves, the largest
group had order $64$.

In this paper we present many explicit families of curves of genus~$2$ and $3$ whose
Jacobians possess large rational torsion subgroups.  The strategy behind our
constructions is to take a product of elliptic curves, each with large 
rational torsion, and to find a curve whose Jacobian is isogenous to
the given product.  Thus it is no surprise that the groups we list
occur as torsion groups of abelian varieties; rather, the point of interest 
is that they occur as torsion groups of Jacobian varieties.

For curves of genus~$2$, we have the following result:

\begin{theorem}
\label{genus2thm}
For every abstract group $G$ listed in the first column of Table~\ref{groups},
there exists a family of curves over $\Q$ of genus~$2$,
parameterized by the rational points on a non-empty Zariski-open subset of a
variety of the type listed in the third column,
whose Jacobians contain a group of rational points isomorphic to $G$.
\end{theorem}

\begin{table}
\begin{center}
\begin{tabular}{|c|c|c|}
$G$			& $|G|$	& Parameterizing variety \\ \hline \hline
$\Z/20\Z$
& 20	& $\PP^2$				\buff\\ \hline
$\Z/21\Z$
& 21	& $\PP^2$				\buff\\ \hline
$\Z/3\Z \times \Z/9\Z$
& 27	& $\PP^2$				\buff\\ \hline
$\Z/30\Z$
& 30	& $\PP^2$				\buff\\ \hline
$\Z/35\Z$
& 35	& positive rank elliptic curve		\buff\\ \hline
$\Z/6\Z \times \Z/6\Z$
& 36	& $\PP^2$				\buff\\ \hline
$\Z/3\Z \times \Z/12\Z$
& 36	& $\PP^2$				\buff\\ \hline
$\Z/40\Z$
& 40	& positive rank elliptic surface	\buff\\ \hline
$\Z/45\Z$
& 45	& positive rank elliptic curve		\buff\\ \hline
$\Z/2\Z \times \Z/24\Z$
& 48	& $\PP^2$				\buff\\ \hline
$\Z/7\Z \times \Z/7\Z$
& 49	& $\PP^0$				\buff\\ \hline
$\Z/5\Z \times \Z/10\Z$
& 50	& positive rank elliptic surface	\buff\\ \hline
$\Z/60\Z$
& 60	& positive rank elliptic curve		\buff\\ \hline
$\Z/63\Z$
& 63	& $\PP^0$				\buff\\ \hline
$\Z/8\Z \times \Z/8\Z$
& 64	& $\PP^2$				\buff\\ \hline
$\Z/2\Z \times \Z/4\Z \times \Z/8\Z$
& 64	& $\PP^2$				\buff\\ \hline
$\Z/6\Z \times \Z/12\Z$
& 72	& positive rank elliptic surface	\buff\\ \hline
$\Z/2\Z \times \Z/6\Z \times \Z/6\Z$
& 72	& positive rank elliptic surface	\buff\\ \hline
$\Z/2\Z \times \Z/2\Z \times \Z/24\Z$
& 96	& positive rank elliptic curve		\buff\\ \hline
$\Z/2\Z \times \Z/2\Z \times \Z/4\Z \times \Z/8\Z$
& 128	& positive rank elliptic surface	\buff\\ \hline
\end{tabular}
\end{center}
\vspace{1ex}
\caption{Families of curves over $\Q$ of genus~$2$ such that $G$ is contained
in the torsion subgroup of the Jacobian.}
\label{groups}
\end{table}

When we say that a family is parameterized by the rational points in a
non-empty Zariski open subset $U$ of a variety $X$, we mean in
particular that the closure of the image of $U$ in the moduli space of
genus-$2$ curves is of the same dimension as $X$.
Also, by ``positive rank elliptic surface'' we mean an elliptic surface
over $\PP^1$ with positive rank.
Note that a family parameterized by the rational points in a non-empty
open subset of $\PP^0$ consists of a single curve.
We will often refer to families parameterized by $\PP^1$ and $\PP^2$
as $1$- and $2$-parameter families.

A similar table expresses our results for curves of genus~$3$.

\begin{theorem}
\label{genus3thm}
For every abstract group $G$ listed in the first column of Table~\ref{groups2},
there exists a family of curves over $\Q$ of genus~$3$,
parameterized by the rational points on a non-empty Zariski-open subset of a
variety of the type listed in the third column,
whose Jacobians contain a group of rational points isomorphic to $G$.
The fourth column of the table indicates whether or not the family consists entirely of
hyperelliptic curves.
\end{theorem}

\begin{table}
\begin{center}
\begin{tabular}{|c|c|c|c|}
$G$			& $|G|$	& Parameterizing variety & All hyp.?\\ \hline \hline
$\Z/2\Z \times \Z/30\Z$
& $60$  & positive rank elliptic curve			& yes\buff\\ \hline
$\Z/10\Z \times \Z/10\Z$
& $100$ & $\PP^1$					& yes\buff\\ \hline
$\Z/2\Z  \times \Z/8\Z  \times \Z/8\Z$
& $128$ & positive rank elliptic surface		& yes\buff\\ \hline
$\Z/4\Z  \times \Z/4\Z  \times \Z/8\Z$
& $128$ & $\PP^1$					& yes\buff\\ \hline
$\Z/4\Z \times \Z/40\Z$	
& $160$ & positive rank elliptic curve			& no\buff\\ \hline
$\Z/2\Z \times \Z/4\Z  \times \Z/24\Z$
& $192$ & positive rank elliptic curve			& no\buff\\ \hline
$\Z/2\Z \times \Z/2\Z  \times \Z/2\Z \times \Z/24\Z$
& $192$ & positive rank elliptic surface		& yes\buff\\ \hline
$\Z/10\Z \times \Z/20\Z$
& $200$ & $\PP^2$					& no\buff\\ \hline
$\Z/6\Z  \times \Z/6\Z  \times \Z/6\Z$
& $216$ & positive rank elliptic curve			& no\buff\\ \hline
$\Z/4\Z \times \Z/60\Z$
& $240$ & positive rank elliptic curve			& no\buff\\ \hline
$\Z/4\Z  \times \Z/8\Z  \times \Z/8\Z$
& $256$ & positive rank elliptic curve			& no\buff\\ \hline
$\Z/2\Z  \times \Z/2\Z  \times \Z/8\Z \times \Z/8\Z$
& $256$ & $\PP^2$					& no\buff\\ \hline
$\Z/2\Z  \times \Z/4\Z  \times \Z/4\Z \times \Z/8\Z$
& $256$ & $\PP^2$					& no\buff\\ \hline
$\Z/2\Z  \times \Z/2\Z  \times \Z/2\Z \times \Z/4\Z \times \Z/8\Z$
& $256$ & $\PP^2$					& yes\buff\\ \hline
$\Z/2\Z \times \Z/12\Z \times \Z/12\Z$
& $288$ & $\PP^2$					& no\buff\\ \hline
$\Z/2\Z \times \Z/2\Z  \times \Z/6\Z \times \Z/12\Z$
& $288$ & positive rank elliptic surface		& yes\buff\\ \hline
$\Z/2\Z  \times \Z/2\Z  \times \Z/4\Z \times \Z/4\Z \times \Z/8\Z$
& $512$ & positive rank elliptic curve			& no\buff\\ \hline
$\Z/2\Z  \times \Z/2\Z  \times \Z/2\Z \times \Z/2\Z \times \Z/4\Z \times \Z/8\Z$
& $512$ & $\PP^1$					& yes\buff\\ \hline
$\Z/6\Z \times \Z/12\Z \times \Z/12\Z$
& $864$ & $\PP^0$					& no\buff\\ \hline
\end{tabular}
\end{center}
\vspace{1ex}
\caption{Families of curves over $\Q$ of genus~$3$ such that $G$ is contained
in the torsion subgroup of the Jacobian.
The final column indicates whether or not
the family consists entirely of hyperelliptic curves. }
\label{groups2}
\end{table}

In Part~\ref{genus1sect} of the paper we review the results on elliptic curves that we will
need to prove these theorems.  In Part~\ref{genus2sect} we show how, given a pair of
non-isomorphic elliptic curves whose Galois modules of $2$-torsion points are isomorphic,
one can construct explicitly a curve of genus~$2$ whose Jacobian is isogenous to
the product of the given elliptic curves.
After giving some quick applications of the construction to the problems of
finding genus-$2$ curves of low conductor and of high rank,
we give a modular interpretation of our construction in Section~\ref{modular}.
The rest of Part~\ref{genus2sect} is taken up with the proof of Theorem~\ref{genus2thm}.
In Part~\ref{genus3sect} we begin with another explicit construction:  We show
in Section~\ref{genus3} how one can construct a curve of genus~$3$ whose Jacobian is isogenous 
to a product of three given elliptic curves, provided that each of the elliptic curves
has a rational $2$-torsion point, that the product of their discriminants is a square,
and that a certain explicitly calculable number depending on the curves is a square.
The remainder of Part~\ref{genus3sect} contains the proof of
Theorem~\ref{genus3thm}.
The reader should note that our verifications of the many entries in 
Tables~\ref{groups} and~\ref{groups2} are organized
not by the sequence of the entries in the tables but rather by the type
of argument the verifications require.  Consequently, the proofs of the theorems
are distributed among several sections.

Throughout the paper, and often without further mention,
we will make use of van Hoeij's Maple package \verb+IntBasis+
for computing Weierstrass models for genus-$1$ curves with a rational point;
Cremona's programs \verb+findinf+ and \verb+mwrank+ for finding points on, and
computing ranks of, elliptic curves over $\Q$;
Mathematica;
and especially PARI.

\section{Genus one}
\label{genus1sect}

In this section we record facts about torsion of elliptic curves over $\Q$ 
that we will need later.
Mazur's theorem~\cite{mazur} states that if $E$ is an elliptic curve over $\Q$, 
then the group of rational torsion points on $E$ is isomorphic to $\Z/N\Z$ with 
$N \le 10$ or $N=12$, or isomorphic to $\Z/2\Z \times \Z/2N\Z$ with $N \le 4$.
For each possibility where the group is not killed by~$2$, the elliptic curves 
having that group as torsion subgroup form a $1$-parameter family.
We will need to have an explicit equation for the universal curve for each family.
For $N=3$, this universal elliptic curve is $y^2=x^3+(x+t)^2/4$ and 
a $3$-torsion point is $(0,t/2)$.
For the other cases, we copy\footnote{
   Actually, we have done a tiny bit more than copy:
   we have expanded the implicit expressions for the parameters $b$ and $c$ 
   in~\cite{kubert} to express $b$ and $c$ in terms of a
   single parameter $t$.}
Table~3 in~\cite{kubert} to our Table~\ref{universal}.

\begin{table}
\begin{center}
\begin{tabular}{|c||c|c|}
$N$ or $(2,2N)$	& $b$	& $c$	\\ \hline \hline
4	& $t$		& $0$		\buff\\ \hline
5	& $t$		& $t$		\buff\\ \hline
6	& $t^2+t$	& $t$		\buff\\ \hline
7	& $t^3-t^2$	& $t^2-t$	\buff\\ \hline
8	& $2t^2-3t+1$	& $\displaystyle \frac{2t^2-3t+1}{t}$	\bbuff\\ \hline
9	& $t^5-2t^4+2t^3-t^2$	& $t^3-t^2$	\buff\\ \hline
10	& $\displaystyle \frac{2t^5-3t^4+t^3}{(t^2-3t+1)^2}$	& $\displaystyle \frac{-2t^3+3t^2-t}{t^2-3t+1}$	\bbuff\\ \hline
12	& $\displaystyle \frac{12t^6-30t^5+34t^4-21t^3+7t^2-t}{(t-1)^4}$	& $\displaystyle \frac{-6t^4+9t^3-5t^2+t}{(t-1)^3}$	\bbuff\\ \hline
(2,4)	& $\displaystyle t^2-\frac{1}{16}$	& $0$		\bbuff\\ \hline
(2,6)	& $\displaystyle \frac{-2t^3+14t^2-22t+10}{(t+3)^2(t-3)^2}$	& $\displaystyle \frac{-2t+10}{(t+3)(t-3)}$		\bbuff\\ \hline
(2,8)	& $\displaystyle \frac{16t^3+16t^2+6t+1}{(8t^2-1)^2}$	& $\displaystyle \frac{16t^3+16t^2+6t+1}{2t(4t+1)(8t^2-1)}$		\bbuff\\ \hline
\end{tabular}
\end{center}
\vspace{1ex}
\caption{Parameters $b,c$ for the universal elliptic curve $y^2+(1-c)xy-by=x^3-bx^2$ over $X_1(N)$ or $X_1(2,2N)$.  In each case, $(0,0)$ is a torsion point of maximal order.}
\label{universal}
\end{table}

Let $E_N^t$ denote the elliptic curve with a rational $N$-torsion point with 
parameter $t$, and similarly define $E_{2,2N}^t$.
We will need to know something about the field of definition of the 
$2$-torsion points on the curves $E_N^t$.
Therefore we record the discriminant $\Delta_N(t)$ of $E_N^t$ modulo squares 
in $\Q(t)$ in Table~\ref{discriminants}.
If $N$ is odd, the discriminant $\Delta_N(t)$ is equal (modulo squares) to 
the discriminant of the cubic field obtained by adjoining the coordinates of 
one $2$-torsion point;
if $N$ is even, $\Delta_N(t)$ is equal (modulo squares) to the discriminant 
of the quadratic field obtained by adjoining the coordinates of a non-rational 
$2$-torsion point.

\begin{table}
\begin{center}
\begin{tabular}{|c||c|}
$N$	& Discriminant $\Delta_N(t)$ modulo squares	\\ \hline \hline
3	& $t(1-27t)$				\buff\\ \hline
4	& $16t+1$				\buff\\ \hline
5	& $t(t^2-11t-1)$			\buff\\ \hline
6	& $(t+1)(9t+1)$				\buff\\ \hline
7	& $t(t-1)(t^3-8t^2+5t+1)$		\buff\\ \hline
8	& $8t^2-8t+1$				\buff\\ \hline
9	& $t(t-1)(t^2-t+1)(t^3-6t^2+3t+1)$	\buff\\ \hline
10	& $(2t-1)(4t^2-2t-1)$			\buff\\ \hline
12	& $(2t^2-2t+1)(6t^2-6t+1)$		\buff\\ \hline
\end{tabular}
\end{center}
\vspace{1ex}
\caption{Discriminant (modulo squares) of the elliptic curve $E_N^t$.}
\label{discriminants}
\end{table}

We will need to know the $x$-coordinates of the nonzero $2$-torsion points
on $E_{2,2N}^t$, at least for $N=3$ and $N=4$.
These are given in Table~\ref{explicit2torsion}.
For $N=4$, the point $T_1$ is the one that is~$4$ times a rational $8$-torsion point.
Note that these $x$-coordinates are also valid for the model
	$$y^2=x^3-bx^2 + [(1-c)x-b]^2/4$$
obtained by completing the square in $y$.

\begin{table}
\begin{center}
\begin{tabular}{|c||c|c|c|}
Curve	& $x(T_1)$	& $x(T_2)$	& $x(T_3)$	\\ \hline \hline
$E_{2,6}^t$	& $\displaystyle \frac{-2t+10}{(t+3)(t-3)}$	& $\displaystyle \frac{-t^3+7t^2-11t+5}{4(t+3)(t-3)^2}$	& $\displaystyle \frac{-2t^2+4t-2}{(t+3)^2(t-3)}$	\bbuff\\ \hline
$E_{2,8}^t$	& $\displaystyle \frac{16t^3+12t^2+2t}{(8t^2-1)^2}$	& $\displaystyle \frac{32t^3+24t^2+8t+1}{16t^2(8t^2-1)}$	& $\displaystyle \frac{-32t^4-32t^3-12t^2-2t}{(4t+1)^2(8t^2-1)}$	\bbuff\\ \hline
\end{tabular}
\end{center}
\vspace{1ex}
\caption{The $x$-coordinates of the $2$-torsion points on $E_{2,2N}^t$.}
\label{explicit2torsion}
\end{table}

For our work with genus-$3$ curves, we will require different
models of the universal elliptic curves $E_N^t$ for $N=4,6,8,10,12$
and $E_{2,2N}^t$ for $N=2,3,4$; in particular, we will want to 
have each curve written in the form $y^2 = x(x^2 + Ax + B)$, where $(x,y)=(0,0)$
is a specified $2$-torsion point.  Table~\ref{universal2} lists the
values of $A$ and $B$ for the curves we will need, as well as the
value of the number $\Delta = A^2-4B$, which differs from the
discriminant of the elliptic curve by a factor of $16B^2$.
The entries in the table were calculated by completing the square in $y$
for the models of the $E_{2N}^t$ and $E_{2,2N}^t$ given above, moving a
rational $2$-torsion point to $x=0$, scaling with respect to $x$ to clear
denominators, and making a linear change of variables in $t$ so as to simplify
the resulting polynomials.  We will refer to these models as $F_{2N}^t$ and
$F_{2,2N}^t$ depending on which universal elliptic curve they model.
However, there are two essentially different ways of putting the curves 
$E_{2,4}^t$ and $E_{2,8}^t$ into the desired form, because one of the $2$-torsion
points on these curves is a multiple of a point of order $4$ while the others
are not.  We  denote the models in which the $2$-torsion point at $x=0$ is a
multiple of a $4$-torsion point by $F_{2,4}^t$ and $F_{2,8}^t$,
and we denote the other models by $F_{4,2}^t$ and $F_{8,2}^t$.

\begin{table}
\begin{center}
\begin{tabular}{|c||l|}
$2N$, $(2,2N)$, or $(2N,2)$	& $A$, $B$, and $\Delta = A^2 -4B$ \\		\hline	\hline
	& $A=2t+1$								\buff\\
4	& $B=t^2$								\buff\\
	& $\Delta = 4t+1$							\buff\\	\hline
	& $A=2t^2+2$								\buff\\
(2,4)	& $B=(t-1)^2(t+1)^2$							\buff\\
	& $\Delta = 16t^2$							\buff\\	 \hline
	& $A=-t^2-6t-1$								\buff\\
(4,2)	& $B=4t(t+1)^2$								\buff\\
	& $\Delta=(t-1)^4$							\buff\\	 \hline
	& $A=-3t^2+6t+1$							\buff\\
6	& $B=-16t^3$								\buff\\
	& $\Delta = (9t+1)(t+1)^3$						\buff\\	 \hline
	& $A=-2t^4 + 12t^2 + 6$							\buff\\
(2,6)	& $B=(t+3)(t-3)(t+1)^3 (t-1)^3$						\buff\\
	& $\Delta = 256 t^2$							\buff\\	\hline
	& $A=2t^4 + 4t^2 - 2$							\buff\\ 
8	& $B=(t+1)^4 (t-1)^4$							\buff\\
	& $\Delta=16(2t^2-1)t^4$						\buff\\	\hline
	& $A=t^8 - 4t^6 + 22t^4 - 4t^2 + 1$					\buff\\
(2,8)	& $B=16 t^4 (t + 1)^4 (t - 1)^4$					\buff\\ 
	& $\Delta = (t^2 - 2t - 1)^2 (t^2 + 2t - 1)^2 (t^2 + 1)^4$		\buff\\	\hline
	& $A = -2t^8 + 8t^6 + 4t^4 + 8t^2 - 2$					\buff\\
(8,2)	& $B = (t^2 - 2t - 1)  (t^2 + 2t - 1) (t^2 + 1)^2 (t + 1)^4 (t - 1)^4$	\buff\\
	& $\Delta = 256 t^8 {}$							\buff\\	\hline
	& $A = -(2t^2 - 2t + 1) (4t^4 - 12t^3 + 6t^2 + 2t - 1)$			\buff\\ 
10	& $B= 16 (t^2 - 3t + 1) (t - 1)^5 t^5$ 					\buff\\ 
	& $\Delta = (4t^2 - 2t - 1) (2t - 1)^5$ 				\buff\\ \hline
	& $A=24t^8 -96t^7 +216t^6 -312t^5 + 288t^4 - 168t^3 + 60t^2 - 12t + 1$	\buff\\
12	& $B=16 (3t^2 - 3t + 1)^2 (t - 1)^6 t^6$				\buff\\
	& $\Delta = (6t^2 - 6t + 1) (2t^2 - 2t + 1)^3 (2t - 1)^6$ 	\buff\\	\hline
\end{tabular}
\end{center}
\vspace{1ex}
\caption{Parameters $A,B$ for the universal elliptic curve
$y^2=x(x^2+Ax+B)$ over $X_1(2N)$ or $X_1(2,2N)$.
The $2$-torsion point $(0,0)$ is twice a rational $4$-torsion point
for the entries marked $(2,4)$ and $(2,8)$, and is not for the entries marked $(4,2)$ and $(8,2)$.}
\label{universal2}
\end{table}

For convenience, we list in Table~\ref{maximal} the coordinates
of a torsion point of maximal order on the curves $F_{2N}^t$, $F_{2,2N}^t$,
and $F_{2N,2}^t$.  Also, in Table~\ref{explicit2torsion2} we give the
$x$-coordinates of the $2$-torsion points other than $(0,0)$ on the 
curves whose $2$-torsion points are all rational.
In the entries for $F_{4,2}^t$ and $F_{8,2}^t$, the point labeled
$T_1$ is twice a rational $4$-torsion point.

\begin{table}
\begin{center}
\begin{tabular}{|c||l|}
$N$, $(2,2N)$, or $(2N,2)$	& $(x,y)$-coordinates \\			\hline	\hline
$4$	& $x=-t$ \buff\\				& $y=t$ \buff\\ \hline
$(2,4)$	& $x=-(t+1)(t-1)$ \buff\\			& $y=2(t+1)(t-1)$ \buff\\			\hline
$(4,2)$	& $x=2(t+1)$ \buff\\				& $y=2(t+1)(t-1)$ \buff\\			\hline
$6$	& $x=-4t$	 \buff\\			& $y=4t(t+1)$ \buff\\				\hline
$(2,6)$	& $x=(t-3)(t+3)(t-1)(t+1)$ \buff\\		& $y=4(t-3)(t+3)(t-1)(t+1)$ \buff\\		\hline
$8$	& $x= -(t+1)^3(t-1)$ \buff\\			& $y=2(t+1)^3 (t-1) t$ \buff\\			\hline
$(2,8)$	& $x=-4(t-1)t(t+1)^3$ \buff\\			& $y=4(t-1)t(t+1)^3(t^2+1)(t^2-2t-1)$ \buff\\	\hline
$(8,2)$	& $x=(t^2+1)(t^2-2t-1)(t-1)(t+1)^3$ \buff\\	& $y=4(t^2+1)(t^2-2t-1)(t-1)(t+1)^3t$ \buff\\	\hline
$10$	& $x=4(t-1)(t^2-3t+1)t^3$ \buff\\		& $y=4(t-1)(t^2-3t+1)t^3(2t-1)$ \buff\\		\hline
$12$	& $x=-4(t-1)(3t^2-3t+1)t^5$ \buff\\		& $y=4(t-1)(3t^2-3t+1)t^5(2t^2-2t+1)(2t-1)$ \buff\\	\hline
\end{tabular}
\end{center}
\vspace{1ex}
\caption{Coordinates of a torsion point of maximal order on 
the universal curves $F_N^t$, $F_{2,2N}^t$, and $F_{2N,2}^t$.}
\label{maximal}
\end{table}

\begin{table}
\begin{center}
\begin{tabular}{|c||c|c|}
$(2,2N)$ or $(2N,2)$	& $x(T_1)$	& $x(T_2)$ \\			\hline	\hline
$(2,4)$	&  $-(t-1)^2$		& $-(t+1)^2$ \buff\\				\hline
$(4,2)$	&  $(t+1)^2$		& $4t$ \buff\\				\hline
$(2,6)$	&  $(t+3)(t-1)^3$	& $(t-3)(t+1)^3$ \buff\\			\hline
$(2,8)$	&  $-16t^4$		& $-(t-1)^4(t+1)^4$ \buff\\			\hline
$(8,2)$	&  $(t-1)^4 (t+1)^4$	& $(t^2+2t-1)(t^2-2t-1)(t^2+1)^2$ \buff\\	\hline
\end{tabular}
\end{center}
\vspace{1ex}
\caption{The $x$-coordinates of the $2$-torsion points on $F_{2,2N}^t$ and $F_{2N,2}^t$
other than $(0,0)$.}
\label{explicit2torsion2}
\end{table}

Finally, we note that while there is no universal elliptic curve over the
modular curve $X(2)$ there is a replacement that will suffice for our purposes.
If $k$ is a field of characteristic different from $2$,
then every elliptic curve over $k$ that has all of
its $2$-torsion defined over $k$ is isomorphic to a twist of 
a specialization of the curve $F_{2,2}^t$ over $k(t)$ defined by 
$y^2 = x(x^2 + Ax + B)$, where $A=-1-t$ and $B=t$.

\section{Genus two}
\label{genus2sect}
\subsection{Conventions}
\label{conventions}

All curves are supposed to be nonsingular and irreducible
unless we specifically mention that they might not be.
The modular curves we consider in section~\ref{modular} are possibly singular.
If $A$ is a variety over a field $k$ and if $K$ is an extension field of $k$,
we will denote by $A_K$ the $K$-scheme $A\times_{\Spec k} \Spec K$.
If $A$ is an abelian variety over a field $k$ and $N$ is a positive integer, 
we will denote by $A[N]$ the $k$-group scheme that is the kernel of the
multiplication-by-$N$ map on $A$.

\subsection{Jacobians $(2,2)$-isogenous to a product of elliptic curves}
\label{construction2}

In this section we will show how one can construct a curve of genus~$2$
whose Jacobian is $(2,2)$-isogenous to a product of two given
elliptic curves, provided one has an isomorphism of their $2$-torsion
groups that does not come from an isomorphism of elliptic curves.
Related results, some of them constructive, have appeared in the
recent literature --- see for example \cite{frey}, \cite{freykani}, \cite{howe},
\cite{kani-number}, \cite{kuhn} --- and in fact the first examples
of curves like the ones we construct go back to Jacobi.
Chapter XXII of~\cite{baker} provides inroads to the older literature
on the subject.  

Suppose $E$ and $F$ are elliptic curves over a separably closed
field $K$, and let $N$ be a positive integer
not divisible by the characteristic of $K$.
The product of the canonical polarizations on $E$ and $F$ is
a principal polarization $\lambda$ on the product variety $A=E\times F$,
and by combining the
Weil pairings on $E[N]$ and $F[N]$ we get a non-degenerate 
alternating pairing $e_N$ from the $N$-torsion of $A$
to the group-scheme of $N$th roots of unity over $K$.
Suppose $G$ is a sub-group-scheme of $A[N]$ that is isotropic with
respect to the pairing $e_N$ and that is maximal with respect to
this property. Then the polarization $N\lambda$ of $A$ reduces
to a principal polarization $\mu$ on the quotient abelian variety
$B=A/G$ (see \cite{milne}, Proposition~16.8, p.~135).
The polarized variety $(B,\mu)$
will be either the polarized Jacobian of a curve over $K$ or
the product of two polarized elliptic curves over $K$.
Suppose $N=2$; in this case it is easy to show that if $(B,\mu)$
is a Jacobian then $G$ must be the graph of an isomorphism
$E[N](K)\rightarrow F[N](K)$. Our first result is that the converse of this
statement is almost true.

\begin{prop}
\label{isotropic}
Let $E$ and $F$ be elliptic curves over a field $k$ whose 
characteristic is not $2$, let $K$ be a separable closure
of $k$, let $A$ be the polarized
abelian surface $E\times F$, and let 
$G\subseteq A[2](K)$ be the graph of a group isomorphism
$\psi\colon E[2](K)\rightarrow F[2](K)$.  Then $G$ is a maximal isotropic
subgroup of $A[2](K)$.  Furthermore, the quotient polarized
abelian variety $A_K / G$ is isomorphic to the polarized Jacobian of 
a curve $C$ over $K$, unless $\psi$ is the restriction to $E[2](K)$
of an isomorphism $E_K\rightarrow F_K$. If $\psi$ gives rise to a curve $C$,
then $C$ and the isomorphism $\Jac C \cong A_K/G$ can be defined
over $k$ if and only if $G$ can be defined over $k$, if and only if
$\psi$ is an isomorphism of Galois modules.
\end{prop}

\begin{proof}
All of the proposition except for the final sentence is the
special case $N=2$ of the results of \cite{kani-number}.  The final statement of
the proposition follows from standard descent arguments that make
use of the fact that the automorphism group of $C$ is naturally 
isomorphic to that of the polarized variety $A_K/G$.
\end{proof}

Let $k$ and $K$ be as in Proposition~\ref{isotropic} and
let $E$ and $F$ be the elliptic curves over $k$ defined
by the equations $y^2 = f$ and $y^2=g$, respectively,
where $f$ and $g$ are separable monic cubic polynomials in $k[x]$
with discriminants $\Delta_f$ and $\Delta_g$. 
Suppose $\psi$ is a Galois-module isomorphism $E[2](K)\rightarrow F[2](K)$
that does not come from an isomorphism $E_K\rightarrow F_K$.
Our next proposition shows how we can use $f$, $g$, and $\psi$ to find
a model for the curve $C$ over $k$ that appears in
Proposition~\ref{isotropic}.

\begin{prop}
\label{explicit}
With notation as above, 
let $\alpha_1$, $\alpha_2$, and $\alpha_3$ be the roots of $f$ in $K$ and
let  $\beta_1$,  $\beta_2$, and  $\beta_3$ be the roots of $g$ in $K$.
Suppose the roots are indexed so that $\psi((\alpha_i,0)) = (\beta_i,0)$. 
The numbers $a_1$, $b_1$, $a_2$, and $b_2$ defined by 
\begin{align*}
a_1 &=   {(\alpha_3-\alpha_2)^2 / ( \beta_3- \beta_2)}
       + {(\alpha_2-\alpha_1)^2 / ( \beta_2- \beta_1)}
       + {(\alpha_1-\alpha_3)^2 / ( \beta_1- \beta_3)} \\
b_1 &=   {( \beta_3- \beta_2)^2 / (\alpha_3-\alpha_2)}
       + {( \beta_2- \beta_1)^2 / (\alpha_2-\alpha_1)}
       + {( \beta_1- \beta_3)^2 / (\alpha_1-\alpha_3)} \\
a_2 &=   \alpha_1( \beta_3- \beta_2) 
       + \alpha_2( \beta_1- \beta_3)
       + \alpha_3( \beta_2- \beta_1) \\
b_2 &=    \beta_1(\alpha_3-\alpha_2)
       +  \beta_2(\alpha_1-\alpha_3)  
       +  \beta_3(\alpha_2-\alpha_1)
\end{align*}
are nonzero, and the ratios $a_1/a_2$ and $b_1/b_2$ are in $k$.
Let $A=\Delta_g a_1/a_2$ and let $B=\Delta_f  b_1/b_2$.
Then the polynomial $h$ defined by 
\begin{align*}
h & = - \bigl(A(\alpha_2-\alpha_1)(\alpha_1-\alpha_3)x^2 + 
              B(\beta_2 -\beta_1) (\beta_1 -\beta_3)\bigr)\\
  & \qquad\cdot \bigl(A(\alpha_3-\alpha_2)(\alpha_2-\alpha_1)x^2 + 
                       B(\beta_3 -\beta_2) (\beta_2 -\beta_1)\bigr)\\
  & \qquad\cdot \bigl(A(\alpha_1-\alpha_3)(\alpha_3-\alpha_2)x^2 + 
                       B(\beta_1 -\beta_3)(\beta_3 -\beta_2)\bigr)
\end{align*}
is a separable sextic in $k[x]$, and
the polarized Jacobian of the curve $C$ over $k$
defined by $y^2 = h$ is isomorphic to the quotient of $E\times F$ by
the graph of $\psi$.
\end{prop}

\begin{proof}
Simple algebra shows that if either $a_1$ or $a_2$ were zero
we would have
$$\beta_3 = \alpha_3 \frac{\beta_2-\beta_1}{\alpha_2-\alpha_1}
             + \frac{ \beta_1\alpha_2 - \beta_2\alpha_1}
                                        {\alpha_2-\alpha_1}.$$
But then the automorphism
$$ \Psi\colon z\mapsto  z \frac{\beta_2-\beta_1}{\alpha_2-\alpha_1}
             + \frac{ \beta_1\alpha_2 - \beta_2\alpha_1}
                                        {\alpha_2-\alpha_1}$$
of $\PP^1_K$ would take $\alpha_i$ to $\beta_i$ for $i=1,2,3$ and
would also take $\infty$ to $\infty$, and this would mean that
$\psi$ came from the isomorphism $E_K\rightarrow F_K$ obtained from $\Psi$,
contrary to our hypotheses.  Therefore $a_1$ and $a_2$ are nonzero.
It is easy to check that the ratio $a_1/a_2$ is fixed by the action
of $S_3$ that permutes the indices of the $\alpha$'s and $\beta$'s.
But the Galois equivariance of the map $\psi$ shows that the action
of $\Gal(K/k)$ on $a_1/a_2$ factors through this action of $S_3$, 
so $a_1/a_2$ is an element of $k$.  A similar argument shows that 
$b_1$ and $b_2$ are nonzero and that $b_1/b_2\in k$.

The group $\Gal(K/k)$ acts on $h$ by permuting its factors, 
so $h$ is an element of $k[x]$.  The coefficient of $x^2$ in each
factor is nonzero, so $h$ is a sextic.  
To show that $h$ is separable it will be enough to show that
the polynomial $\gb\in k[u]$ defined by
\begin{align*}
\gb & = - \bigl(A(\alpha_2-\alpha_1)(\alpha_1-\alpha_3)u + 
                       B(\beta_2 -\beta_1) (\beta_1 -\beta_3)\bigr)\\
    & \qquad\cdot \bigl(A(\alpha_3-\alpha_2)(\alpha_2-\alpha_1)u + 
                       B(\beta_3 -\beta_2) (\beta_2 -\beta_1)\bigr)\\
    & \qquad\cdot \bigl(A(\alpha_1-\alpha_3)(\alpha_3-\alpha_2)u + 
                       B (\beta_1 -\beta_3)(\beta_3 -\beta_2)\bigr)
\end{align*}
is separable, because $h(x)=\gb(x^2)$ and the roots of $\gb$ are
nonzero.
Let $t_1=-(A/B)(b_2/b_1)$ and let $t_2$ be the element
$$t_2 = \frac{1}{b_1}
      \left(  \frac{\beta_1(\beta_3-\beta_2)^2}{\alpha_3-\alpha_2}
            + \frac{\beta_2(\beta_1-\beta_3)^2}{\alpha_1-\alpha_3} 
            + \frac{\beta_3(\beta_2-\beta_1)^2}{\alpha_2-\alpha_1}\right)
$$
of $k$.  The reader may verify that the automorphism 
$z\mapsto t_1z+t_2$ of $\PP^1_K$ takes the
roots of $\gb$ to the roots of $g$.  The roots of $g$ are distinct by
assumption, so the roots of $\gb$ must also be distinct, so $\gb$
is separable.

Now we turn to the final statement of the proposition.
Let $\Fb$ be the elliptic curve over $k$ defined by $v^2=\gb$.
Once one knows that $z\mapsto t_1z+t_2$ takes the roots of $\gb$ to
those of $g$, it is a simple matter to verify that the map
$$(u,v) \mapsto (t_1u+t_2, (\Delta_f/B^3) v)$$
provides an isomorphism between $\Fb$ and $F$.
Thus we can define a morphism $\chi\colon C\rightarrow F$ of degree~$2$
by $$(x,y)\mapsto(t_1x^2+t_2, (\Delta_f/B^3) y).$$  The involution $\tau$
of $C$ defined by this double cover is given by $(x,y)\mapsto(-x,y)$.

Let $\Eb$ be the elliptic curve over $k$ defined
by $v^2 = \fb$, where $\fb\in k[u]$ is given by 
\begin{align*}
\fb & = - \bigl(A(\alpha_2-\alpha_1)(\alpha_1-\alpha_3) + 
                       B(\beta_2 -\beta_1) (\beta_1 -\beta_3)u\bigr)\\
    & \qquad\cdot \bigl(A(\alpha_3-\alpha_2)(\alpha_2-\alpha_1) + 
                       B(\beta_3 -\beta_2) (\beta_2 -\beta_1)u\bigr)\\
    & \qquad\cdot \bigl(A(\alpha_1-\alpha_3)(\alpha_3-\alpha_2) + 
                       B (\beta_1 -\beta_3)(\beta_3 -\beta_2)u\bigr).
\end{align*}
The $\alpha$--$\beta$ symmetry in our equations shows that there
is an isomorphism $\Eb\rightarrow E$ given by
$$(u,v) \mapsto (s_1u+s_2, (\Delta_g/A^3) v),$$
where $s_1$ and $s_2$ are the elements of $k$ defined by exchanging
$\alpha$'s and $\beta$'s in the definitions of $t_1$ and $t_2$.
Thus we get a $k$-morphism $\varphi\colon C\rightarrow E$ of degree $2$
defined by $$(x,y)\mapsto(s_1/x^2+s_2, (\Delta_g/A^3) (y/x^3)).$$
The involution $\sigma$ of $C$ defined by this double cover is given
by $(x,y)\mapsto(-x,-y)$.

Let $A=E\times F$, let $J$ be the Jacobian of $C$, and let 
$\omega \colon A\rightarrow J$ be the morphism $\varphi^*\times\chi^*$.
Note that the image of $\varphi^*$ in $J$ is fixed by $\sigma^*$, while
the image of $\chi^*$ is fixed by $\tau^*$; since $\sigma^*\tau^*=-1$,
we see that $\omega$ is an isogeny.
Let $\mu$ be the canonical polarization of $J$.
The fact that $\varphi$ has degree $2$ implies that
$\widehat{\varphi^*}\mu\varphi^*$ is the multiplication-by-$2$ map on $E$,
and similarly $\widehat{\chi^*}\mu\chi^*$ is the multiplication-by-$2$ map 
on $F$; here $\widehat{\ }$ indicates the dual morphism.
If we let $\lambda$ be the product polarization on 
$A$ obtained from the canonical polarizations on $E$ and $F$,
then we have a commutative diagram
$$
\begin{matrix}
A                 & \mapright{2\lambda} & \widehat{A} \\
\Bmapdown{\omega} &                     &\Bmapup{\widehat{\omega}}\\
J                 & \mapright{\mu}      & \widehat{J}.
\end{matrix}
$$
The diagram shows that $\omega$ must have degree four, and its
kernel lies in the $2$-torsion of $A$.
By using the explicit representation of $2$-torsion elements 
of $E$, $F$, and $J$ as degree-zero $K$-divisors on $E$, $F$, and $C$ that
are supported only on Weierstrass points, one may check easily
that the graph $G$ of $\psi$ is contained in $\ker \omega$,
and since $\#G=\#\ker\omega$, we must have $G=\ker\omega$.
\end{proof}

Below we give a few quick applications of Proposition~\ref{explicit}.
First, we exhibit a curve of genus~$2$ over $\Q$ whose Jacobian has a very
small conductor.  Mestre \cite{mestreconducteurs} proved
under standard conjectures
that the conductor of a $g$-dimensional abelian variety over $\Q$
must be greater than $(10.32)^g$, so for a $2$-dimensional variety
a conductor of $121$ is close to the minimum of $107$.

\begin{cor}
\label{conductor121}
The conductor of the Jacobian of the curve $y^2 = -2x^6-10x^4+26x^2+242$
is~$121$.
\end{cor}

\begin{proof}
Take $E$ and $F$ to be the modular curves $X_1(11)$ and $X_0(11)$ over $\Q$.
The $\Q$-rational $5$-isogeny $E\rightarrow F$ gives us
a Galois-module isomorphism 
$\psi\colon E[2](\Qbar)\rightarrow F[2](\Qbar)$, and $\psi$ does not come from
an isomorphism $E_\Qbar\rightarrow F_\Qbar$ because
$E_\Qbar$ and $F_\Qbar$ are not
isomorphic to one another.  Applying Proposition~\ref{explicit} to convenient
models of $E$ and $F$ and simplifying the resulting equation gives us
the curve in the statement of the corollary.
\end{proof}

\begin{rem}
The curve in Corollary~\ref{conductor121} is none other than $X_0(22)$.
An isomorphism from the model
	$$Y^2 = (X^3 + 2X^2 - 4X + 8)(X^3 - 2X^2 + 4X - 4)$$
for $X_0(22)$ given in~\cite{hasegawa} to our curve
	$$y^2 = -2x^6-10x^4+26x^2+242$$
is given by $(x,y)=(1-4/X,16Y/X^3)$.
\end{rem}

Our next two results give sufficient conditions for there to be a genus-$2$
curve whose Jacobian is isogenous to the product of a given elliptic curve
with itself. (The sufficient conditions given by Kani~\cite{kani-existence} apply
only when $k$ is algebraically-closed.)

\begin{cor}
\label{EtimesEcharp}
Suppose $E$ is an elliptic curve over a field $k$ of characteristic $p>2$.
If $j(E)$ is neither $0$ nor $1728$,
then there exists a genus-$2$ curve over $k$ whose Jacobian is
isogenous over $k$ to $E \times E$.
\end{cor}

\begin{proof}
First suppose that $j(E)\not\in\F_p$.
Let $y^2 = f$ be a Weierstrass equation for $E$,
and let $F$ be the elliptic curve over $k$ defined by $y^2 = f^{(p)}$,
where $f^{(p)}$ is the polynomial obtained from $f$ by replacing 
every coefficient with its $p$th power.
Note that $j(F) = j(E)^p \neq j(E)$ since $j(E)\not\in\F_p$.
The Frobenius morphism $(x,y)\mapsto (x^p,y^p)$ from $E$ to $F$
defines a Galois-module isomorphism $\psi\colon E[2](K)\ra F[2](K)$.
The curves $E$ and $F$ have different $j$-invariants and so
are not isomorphic to one another, so Proposition~\ref{isotropic}
provides us with the required genus-$2$ curve.

Next suppose that $j(E)\in \F_p$. 
Suppose we could find a Galois-module automorphism $\psi$ of $E[2](K)$ other than the
identity. Such an automorphism could not come from an automorphism of $E_K$
because our assumption on $j(E)$ implies that $\Aut E_K = \{\pm 1\}$,
so Proposition~\ref{isotropic} would give us a curve whose
Jacobian is $(2,2)$-isogenous to $E\times E$.  Thus we will be done
if we can find such a $\psi$.

If $\#E[2](k)=4$ then $E[2](K)$ is a trivial Galois module, and 
we may take $\psi$ to be any one of the $5$ non-trivial automorphisms of $E[2](K)$.

If $\#E[2](k)=2$, we may take $\psi$ to be the automorphism of $E[2](K)$
that swaps the two non-rational points.

Finally, suppose  $\#E[2](k)=1$. The fact that $\Aut E_K = \{\pm 1\}$ implies
that $E$ has a model of the form $dy^2 = f$, where $f\in \F_p[x]$ and $d\in k$.
The $2$-torsion points of $E$ are defined over the splitting field of $f$,
which in this case is a cyclic cubic extension of $k$ because $f\in \F_p[x]$.
Therefore the two group automorphisms of $E[2](K)$ that rotate the nonzero 
points are actually Galois-module isomorphisms, and we have two choices
for $\psi$.
\end{proof}

In characteristic~$0$ it is a little more difficult to find such curves $C$.

\begin{cor}
\label{etimese}
Let $E$ be an elliptic curve over a field $k$ of characteristic not~$2$,
and suppose $\End E_K=\Z$.
If $E$ contains a $k$-rational finite subgroup-scheme $G$ that is not the kernel of 
multiplication-by-$N$ for any $N$, 
then there exists a genus-$2$ curve over $k$ whose Jacobian is
isogenous over $k$ to $E \times E$.
\end{cor}

\begin{proof}
We may assume that the order of $G$ is a prime $q$.
If $q=2$, then there is a nonzero $k$-rational {\em point} of order~$2$,
and we may use in Proposition~\ref{isotropic} the isomorphism
$\psi: E[2](K) \rightarrow E[2](K)$ that interchanges the
other $2$-torsion points.
If $q$ is an odd prime, then the natural isogeny from $E$ to the elliptic
curve $F=E/G$ over $k$ defines an isomorphism of Galois-modules
$\psi: E[2](K) \rightarrow F[2](K)$.
The condition $\End E_K=\Z$ ensures that $E_K$ and $F_K$ are not isomorphic,
so the result again follows from Proposition~\ref{isotropic}.
\end{proof}

\begin{rem}
The conclusion of Corollary~\ref{etimese} holds for some
elliptic curves in characteristic~$0$ that do not satisfy the
corollary's hypotheses.
For example, let $E$ be the elliptic curve
	$$E: y^2 = x^3-169x+845$$
over $\Q$.
The cubic on the right is irreducible, and has square discriminant $13^4$,
so its Galois group is $A_3$.
Therefore any isomorphism $\psi: E[2](\Qbar) \rightarrow E[2](\Qbar)$
that rotates the three non-trivial $2$-torsion points will be defined
over $\Q$.
Since $j(E) \not=0$, such a rotation cannot be the restriction of an
automorphism of $E$, so by Proposition~\ref{isotropic}, we obtain
a genus-$2$ curve over $\Q$ whose Jacobian is $(2,2)$-isogenous
over $\Q$ to $E \times E$.
On the other hand, $E$ is curve 676D1 in~\cite{cremona},
which has no $\Q$-rational cyclic subgroups.
\end{rem}

We can use Corollary~\ref{etimese} to construct genus-$2$ curves over $\Q$
whose Jacobians have high rank, as was also noticed by St\'efane Fermigier.

\begin{cor}
The Jacobian of the curve
\begin{equation*}
\begin{split}
	y^2	& = -1707131824107329945 \cdot (x^2+55871769054504519799033274614104129) \\
		& \quad \cdot (x^4-1086862437115841494920959046499163042 x^2 \\
		& \quad +3121654577279888882305769763628790308995888274656243920700573254848641)
\end{split}
\end{equation*}
has rank $28$ over $\Q$.
\end{cor}

\begin{proof}
According to~\cite{fermigier}, the elliptic curve
  $$E: y^2 = x(x^2 + 2429469980725060 x + 275130703388172136833647756388)$$
has rank~$14$, and $(0,0)$ is a rational $2$-torsion point on $E$.
The $j$-invariant is
	$$\frac{483941743120924000812123996730853715578647268051688786879688}
           {5250870830712351132421548861849566889806152906127048721},$$
which is not an integer, so $E$ cannot have complex multiplication.
Using Corollary~\ref{etimese} and the formulas of Proposition~\ref{explicit}, 
we obtain the desired genus-$2$ curve over $\Q$ whose Jacobian
is $(2,2)$-isogenous to $E \times E$.
\end{proof}

\subsection{A modular interpretation}
\label{modular}
One of our goals in this paper is to construct curves over $\Q$ of genus~$2$
whose Jacobians have large rational torsion subgroups, and our strategy
will be to use Proposition~\ref{explicit} to ``tie together'' two elliptic
curves that each have large torsion subgroups.  In particular, every
curve $C$ we construct will come equipped with a $(2,2)$-isogeny
$E\times F\rightarrow \Jac C$,
where $E$ and $F$ have some particular rational torsion
structure.  We would like to construct the moduli space of
curves equipped with such isogenies.

Suppose $k$ is a field and $K$ is a separable closure of $k$.
Pick a set of elements $\{\zeta_M : M\in\Z_{>0}\}$ of $K$ such that
$\zeta_M$ generates the group of $M$th roots of unity in $K$
and $\zeta_M=\zeta_{jM}^j$ for all integers $M,j>0$.
By a {\it full level-$M$ structure\/} on an elliptic
curve $E/K$ we mean a pair of points $(P,Q)$ in $E(K)$ 
that form a Drinfeld basis for $E[M]$ (see \cite{katzmazur}, Chapter 1)
and such that
$P$ and $Q$ pair to $\zeta_M$ under the Weil pairing on $E[M]$.
This corresponds to the moduli problem denoted in~\cite{katzmazur}
by $[\Gamma(M)]^\can$ (see~\cite{katzmazur}, Sections~3.1 and~9.1),
but only because we are working over a field
--- we would have to be more careful with the roots of unity 
otherwise.  There is an obvious right action of the
group $\SL_2(\Z/M\Z)$ on the set of
full level-$M$ structures on a given curve $E$.
Suppose $G$ is a subgroup of $\SL_2(\Z/M\Z)$; by a
{\it partial level-$M$ structure of type $G$\/} on a curve $E/K$
we mean a $G$-orbit of full level-$M$ structures on $E$.
If $N$ is a positive divisor of $M$,
then an {\it $(N,M)$-structure\/} on an elliptic curve $E/K$ is a
pair $(P,Q)$ of points on $E(K)$ such that $Q$ has ``exact order $M$''
(see~\cite{katzmazur}, Chapter 1)
and such that $P$ and $(M/N)Q$ form a full level-$N$ structure on $E$;
this is an example of a partial level-$M$ structure.
If $E$ is an elliptic curve over~$k$, then by a {\it partial
level-$M$ structure of type $G$\/} on $E$ we mean a 
partial level-$M$ structure of type $G$ on $E_K$ that is stable
under the action of $\Gal(K/k)$.

We let $X(M)$ denote the usual
compactified coarse moduli space of elliptic curves 
with full level-$M$ structure; we view $X(M)$ as a curve 
over $k(\zeta_M)$. Note that if $\charact k$ divides $M$ then
$X(M)$ will have several components.  For every 
subgroup $G$ of $\SL_2(\Z/M\Z)$ there is also a
modular curve, which we will denote by 
$X(M;G)$, that parameterizes elliptic curves with 
partial level-$M$ structure of type $G$.  The curve $X(M;G)$
is a $k(\zeta_M^a)$-scheme,
where $a\in(\Z/M\Z)^*$ is a generator of the 
subgroup $\det G\subset(\Z/M\Z)^*$.
Finally, we denote by $X_1(N,M)$ the modular curve that parameterizes
elliptic curves with $(N,M)$-structure. The curve $X_1(N,M)$ is a scheme
over $k(\zeta_N)$.

Suppose the characteristic of the base field $k$ is not $2$, and
suppose we are given two integers $M$ and $N$ and subgroups
$G\subset\SL_2(\Z/M\Z)$ and $H\subset\SL_2(\Z/N\Z)$.
Let $\ell$ be the smallest field containing the fields of definition
of $X(M;G)$ and $X(N;H)$.
We are interested in the functor $\CF$ from the category
of fields over $\ell$ to the category of sets defined as follows:
If $r\supset \ell$  is a field with separable closure $R$,
then $\CF(r)$ is the set of all $R$-isomorphism classes of triples
$((E,\alpha),(F,\beta),\psi)$, 
where $E$ is an elliptic curve with 
partial level-$M$ structure $\alpha$ of type $G$ over $r$,
where $F$ is an elliptic curve with 
partial level-$N$ structure $\beta$ of type $H$ over $r$, and where
$\psi$ is a Galois-module isomorphism $E[2](R)\rightarrow F[2](R)$;
here we say that $((E,\alpha),(F,\beta),\psi)$ and 
$((E',\alpha'),(F',\beta'),\psi')$ are $R$-isomorphic
if there are isomorphisms $\varphi\colon (E,\alpha)_R\rightarrow (E',\alpha')_R$ and 
$\chi\colon (F,\beta)_R\rightarrow (F',\beta')_R$ such that 
$\psi'\circ \varphi = \chi\circ\psi$ on $E[2](R)$.
We will show that this functor is represented by the $\ell$-scheme
$Y^0$ defined in the next paragraph.

The modular curve $X(2)$ is defined over $k$, 
and since $\charact k\neq2$ it has only one component.
The covering $X(2)\rightarrow X(1)$ is Galois with group $S=\SL_2(\Z/2\Z)$,
and the action of an element $s\in S$ on a point $X(2)$ is determined
by its action on the triple $(E,P,Q)$ corresponding to that point.
Let $Z_1=X(2)_\ell\times_{X(1)} X(M;G)_\ell$ and let 
$Z_2= X(2)_\ell\times_{X(1)} X(N;H)_\ell$, 
where $\times_{X(1)}$ means the fiber product over $X(1)_\ell$.
The covers  $Z_1\rightarrow X(M;G)_\ell$ and $Z_2\rightarrow X(N;H)_\ell$ are
Galois with group $S$.
Let $Z$ be the $2$-dimensional $\ell$-scheme $Z_1\times Z_2$,
where $\times$ means the fiber product over $\Spec \ell$,
let $S$ act on the cover $Z\rightarrow X(M;G)_\ell\times X(N;H)_\ell$ 
diagonally, and let $Y$ be the quotient surface of $Z$ by this action.
Finally, let $Y^0$ be the open subvariety of $Y$ that lies
over the open subvariety of $X(1)_\ell \times X(1)_\ell$ where
neither factor is $\infty$.

\begin{prop}
\label{moduli}
The scheme $Y^0$ represents $\CF$.
\end{prop}

\begin{proof}
First let us determine how to describe the $r$-valued points on $Y^0$.
Let $y$ be such a point.  Then $y$ corresponds to a $\Gal(R/r)$-stable
$S$-orbit of elements of $Z(R)$ that lie over the finite part of
$X(1)_R\times X(1)_R$.  Let $z$ be one of the points in this orbit.
According to the modular interpretation of $Z$, the point $z$ 
corresponds to the $R$-isomorphism class of a quadruple
$$\big( (E,P,Q), (E,\alpha), (F,U,V),(F,\beta)\big),$$
where $E$ and $F$ are elliptic curves over $R$,
where $P$ and $Q$ are independent $2$-torsion points in $E(R)$ and
$U$ and $V$ are independent $2$-torsion points in $F(R)$,
where $\alpha$ is a partial level-$M$ structure of type $G$ on $E$,
and where $\beta$ is a partial level-$N$ structure of type $H$ on $F$.
When an element of $S$ is applied to this quadruple, the
only things that get changed are the points $P$, $Q$, $U$, and $V$, so
the fact that the $S$-orbit containing $z$ is defined over $r$ means
that $(E,\alpha)$ is isomorphic to all of its Galois conjugates 
and $(F,\beta)$ is isomorphic to all of its Galois conjugates.
According to Proposition~3.2 (p.~274) of ~\cite{delignerapoport},
this means that $(E,\alpha)$ and $(F,\beta)$ can be defined over~$r$.
If we think of $E$ and $F$ as curves over $r$, then the fact that 
the $S$-orbit of $z$ is defined over $r$ means exactly that
the group isomorphism $\psi\colon E[2](R)\rightarrow F[2](R)$ defined by sending
$P$ to $U$ and $Q$ to $V$ is Galois equivariant.  Thus, 
a point $y\in Y^0(r)$ gives us a triple
$((E,\alpha),(F,\beta),\psi)$ --- but only up to $R$-isomorphism.
And clearly the $R$-isomorphism class of such a triple will give us
a point on $Y^0$.  This gives us a bijection between 
$\CF(r)$ and $\Hom(\Spec r,Y^0)$ for every $r$, and
this collection of bijections is easily seen to 
provide a natural equivalence $\CF\leftrightarrow\Hom(\ \cdot\ ,Y^0)$.
\end{proof}

Let $W$ be the open subscheme of $Y^0$ whose $r$-valued points
correspond to $R$-isomorphism classes of triples 
$((E,\alpha),(F,\beta),\psi)$ such that $\psi$ does not come from
an isomorphism between $E_R$ and $F_R$.  From 
Proposition~\ref{isotropic} we 
see that the $r$-valued points of $W$ correspond to $R$-isomorphism
classes of triples $(C, (E,\alpha), (F,\beta))$, where
$(E,\alpha)$ is an elliptic curve with partial level-$M$
structure of type $G$ over $r$,
where $(F,\beta)$ is an elliptic curve with partial level-$N$
structure of type $H$ over $r$, and where
$C$ is a curve of genus~$2$ over $r$ provided with a $(2,2)$-isogeny
$E\times F\rightarrow \Jac C$ that takes twice the
canonical polarization of $E\times F$ to the canonical polarization
of $\Jac C$.  We abbreviate this by saying that $W$ is the moduli
space for such triples.

\begin{cor}
\label{genus2moduli}
Let $M'$ be the least common multiple of $2$ and $M$ and 
let $N'$ be the least common multiple of $2$ and $N$.
Every geometric component of $W$ is an open subvariety
of a quotient surface of $X(M')_K\times X(N')_K$.
\end{cor}

\begin{proof}
Let $Z_1$ and $Z_2$ be as in the construction of $Y$ above
and let $\varphi$ and $\psi$ be the natural quotient maps from $X(M')_K$ to
$X(2)_K$ and to $X(M;G)_K$, respectively.  For every $s$ in the covering
group $S$ of $X(2)/X(1)$ we get a morphism 
$\Phi_s \colon X(M')_K\rightarrow (Z_1)_K$ by combining the morphisms 
$s\varphi$ and $\psi$.
It is clear from the modular interpretation of these schemes
that the maps $\Phi_s$ provide a surjective morphism from
the sum of six copies of $X(M')_K$ to $(Z_1)_K$.
Similarly we find a surjective morphism from
the sum of six copies of $X(N')_K$ to $(Z_2)_K$.
Therefore every component of $Z_K = Z_{1K}\times Z_{2K}$ is a quotient
surface of $X(M')_K\times X(N')_K$, and every component of $W_K$
is an open subvariety of a quotient surface of $X(M')_K\times X(N')_K$.
\end{proof}

We will be interested in finding genus-$2$ curves over $\Q$
whose Jacobians are equipped with $(2,2)$-isogenies from a product
of elliptic curves with specified rational torsion structures.
Thus we will want to look at the $\Q$-rational points on the moduli space $W$,
and it would be particularly nice to find subvarieties of $W$ whose
$\Q$-rational points are Zariski dense.  In the next few sections we will
find such subvarieties for several different choices of torsion structures,
although we will not phrase our arguments in terms of moduli spaces.

\begin{example}
Suppose we are interested in tying together
an elliptic curve with $(2,M)$-structure and an elliptic curve
with $(2,N)$-structure, where $M$ and $N$ are even integers.
It is easy to check that $X(2)\times_{X(1)} X_1(2,M)$ is the
sum of six copies of $X_1(2,M)$, and 
$X(2)\times_{X(1)} X_1(2,N)$ is the sum of six copies of $X_1(2,N)$.
The group $S$ acts on each of these varieties by permuting the summands,
and the quotient surface $Y$ is the sum of six copies of 
$X_1(2,M)\times X_1(2,N)$.
Thus in this case $W$ has a very simple structure.
The reader is encouraged to work out the structure of $W$ for
other pairs of partial level structures and to 
keep the results in mind when reading the following sections.
\end{example}

\subsection{Building Jacobians from elliptic curves $E$ with $\#E(\Q)[2]=4$}

If two elliptic curves over $\Q$ are to have $2$-torsion subgroups
isomorphic as Galois-modules,
it is necessary that they have the same number of rational $2$-torsion points!
In this section we consider the case where this number is~$4$,
so we are concerned with the families of elliptic curves
over $\Q$ with torsion subgroup containing $\Z/2\Z \times \Z/2N\Z$, for $N \le 4$.
Any member of the family with $N=N_1$ can be paired with any
member of the family with $N=N_2$, since the $2$-torsion
subgroups are automatically isomorphic as (trivial) Galois-modules.
Moreover the generic members of each family (choosing a different
indeterminate parameter for each) are clearly not isomorphic to
each other, so by Proposition~\ref{isotropic}, we get $2$-parameter
families of genus-$2$ curves whose Jacobians map via a
$(2,2)$-isogeny to $E_1 \times E_2$.
In other words, we have shown that the corresponding moduli
space is a union of rational surfaces over $\Q$.
(This also follows immediately from the example at the end of the
preceding section.)
That these families really have two parameters can be seen from the
fact that over $\C$, one can specify the $j$-invariants of
the two elliptic curves arbitrarily and independently\footnote{Actually,
one should choose the $j$-invariants to be different, so that the
elliptic curves are guaranteed not to be isomorphic, but
this is an open condition, so the number of parameters is not
reduced by this constraint.}.  Similar arguments apply later in this paper;
we leave the details to the reader.

The product of the rational torsion in the two elliptic curves
does not map injectively to the rational torsion points of the Jacobian,
but only a $(2,2)$-subgroup is killed.
The group structure of the image of this product in the Jacobian
depends on $N_1$ and $N_2$, but also on the choice of isomorphism
between the $2$-torsion of the two curves if $N_1$ and $N_2$ are
even, since in this case each elliptic curve has a $2$-torsion
point which is distinguished by the property of being $N_i$
times another rational torsion point.
For instance, if $N_1=N_2=4$, elementary calculations with abelian
groups show that this group has structure
$\Z/2\Z \times \Z/4\Z \times \Z/8\Z$ if these special $2$-torsion points
are identified under $\psi$, and $\Z/8\Z \times \Z/8\Z$ otherwise.
If $N_1=3$ and $N_2=4$, then we obtain $\Z/2\Z \times \Z/24\Z$.
If $N_1=N_2=3$, then we obtain $\Z/6\Z \times \Z/6\Z$.
We have not considered the cases where $N_i \le 2$,
since these cases lead to subgroups of the above.

Although some rational $2$-power torsion is lost upon passing from
$E \times F$ to the Jacobian, there is also the possibility
that some $2$-power torsion can be gained: a non-rational point
on $E \times F$ might map to a rational point on $J$.
This phenomenon will be explored in Section~\ref{gaining}.

\subsection{Building Jacobians from elliptic curves $E$ with $\#E(\Q)[2]=2$}

We now consider elliptic curves $E$ and $F$ having torsion
subgroups $\Z/N\Z$ and $\Z/N'\Z$ with even $N, N' \le 12$.
An isomorphism of Galois-modules from $E[2](\Qbar)$ to $F[2](\Qbar)$ must
map the rational $2$-torsion point to the rational $2$-torsion point,
so we see that such an isomorphism
exists if and only if the quadratic field over which the
non-rational $2$-torsion points of $E$ are defined equals the
quadratic field for $F$, and this holds if and only if the
discriminants of $E$ and $F$ are equal modulo squares.
We are thus led to the problem of finding the rational solutions to
\begin{equation}
\label{evensurface}
	\Delta_N(t) y^2 = \Delta_{N'}(u)
\end{equation}
outside the $1$-dimensional closed subset corresponding to cases
where $E$ or $F$ degenerates or where $j(E)=j(F)$.
Each such solution gives rise to a Jacobian of a genus-$2$ curve over $\Q$
whose torsion subgroup contains the quotient of $\Z/N\Z \times \Z/N'\Z$
by the identification of the points of order~$2$ in each factor.

If $N'=4$, then~(\ref{evensurface}) is
	$$\Delta_N(t) y^2 = 16u+1,$$
which is a rational surface over $\Q$, since we can solve for $u$ in
terms of $t$ and $y$.
In particular, for $N=10$, we obtain a $2$-parameter family of Jacobians
whose torsion subgroups contain $\Z/20\Z$.
(The other $N$ will give results which are subsumed in our other results.)

If $N'=6$, then~(\ref{evensurface}) is
\begin{equation}
\label{thirty}
	\Delta_N(t) y^2 = (u+1)(9u+1).
\end{equation}
This can be considered as a conic over $\Q(t)$ with a
$\Q(t)$-rational point, namely $(u,y)=(-1,0)$, and it
is easy to see that this makes~(\ref{thirty}) a rational surface.
In particular, for $N=10$ or $N=12$, we obtain a $2$-parameter family
of Jacobians whose torsion subgroups contain
$\Z/30\Z$ or $\Z/3\Z \times \Z/12\Z$, respectively.

If $N'=10$ and $N=8$, then~(\ref{evensurface}) is
	$$(8t^2-8t+1)y^2=(2u-1)(4u^2-2u-1).$$
If we set $t=(s^2-2s+3)/(4s^2+4)$, we obtain a split elliptic surface
over the $s$-line, and the $\Q(s)$-rational point 
$(u,y)=(-1/2,(2s^2+2)/(s^2-2s-1))$ 
is of infinite order, since its
specialization at $s=0$ is of infinite order on the resulting
elliptic curve over $\Q$.
Thus we have an elliptic surface over $\PP^1$ of positive rank,
and the $\Q$-rational points on this surface outside of the
$1$-dimensional set of degenerate solutions parameterize Jacobians
over $\Q$ whose torsion subgroups contain $\Z/40\Z$.

If $N'=N=10$, then~(\ref{evensurface}) is
	$$(2t-1)(4t^2-2t-1) y^2 = (2u-1)(4u^2-2u-1)$$
which is an elliptic surface over the $t$-line, and $(u,y)=(t,1)$
is a $\Q(t)$-rational point of infinite order, since under the obvious
isomorphism over $\Q(t)(\sqrt{(2t-1)(4t^2-2t-1)})$ to the elliptic curve
	$$y^2 = (2u-1)(4u^2-2u-1)$$
it maps to a point with non-constant $u$-coordinate.
Hence we obtain a positive rank elliptic surface whose points
(outside a $1$-dimensional set) parameterize Jacobians
whose torsion subgroups contain $\Z/5\Z \times \Z/10\Z$.

Similarly, if $N'=N=12$, then~(\ref{evensurface}) is
	$$(2t^2-2t+1)(6t^2-6t+1) y^2 = (2u^2-2u+1)(6u^2-6u+1),$$
which again is an elliptic surface over the $t$-line if we
choose $(u,y)=(t,1)$ as the zero section.
We then have the $\Q(t)$-rational point $(u,y)=(t,-1)$, which
is of infinite order, for the same reason as in the previous case.
Hence we obtain a family of Jacobians, parameterized by 
the points on an open subset of a positive rank elliptic surface,
whose torsion subgroups contain $\Z/6\Z \times \Z/12\Z$.

Finally, if $N'=10$ and $N=12$, then~(\ref{evensurface}) is
\begin{equation}
\label{tenandtwelve}
	(2t^2-2t+1)(6t^2-6t+1) y^2 = (2u-1)(4u^2-2u-1).
\end{equation}
If we choose $t=1/3$, the resulting elliptic curve is
curve 900A1 of~\cite{cremona}, which has rank~$1$.
There are only finitely many rational points on this elliptic curve
that give $u$ such that $E_{10}^u$ degenerates or is isomorphic to
$E_{12}^{1/3}$, so we obtain a family of Jacobians, parameterized
by the points on an open subset of a positive rank
elliptic curve, whose torsion subgroups contain $\Z/60\Z$.

\begin{rem}
In fact, there are infinitely many other specializations of $t$ for
which~(\ref{tenandtwelve}) becomes an elliptic curve of positive rank.
\end{rem}

\subsection{Building Jacobians from elliptic curves $E$ with $\#E(\Q)[2]=1$}

Here we consider elliptic curves $E$ and $F$ having torsion
subgroups $\Z/N\Z$ and $\Z/N'\Z$, respectively, with $N,N'$ odd
(and at most~$9$).
For an elliptic curve $y^2=f(x)$ with trivial rational $2$-torsion,
each non-trivial $2$-torsion point is defined over a cubic extension,
namely the extension obtained by adjoining a root of $f(x)$.
The Galois-modules $E[2](\Qbar)$ and $F[2](\Qbar)$ are isomorphic
if and only if the corresponding cubic fields are isomorphic.
In this case, the discriminants of the elliptic curves must
be equal modulo squares.
The converse is not quite true (cubic fields having the same
discriminant modulo squares are not necessarily isomorphic),
but it will turn out that the discriminants often contain
enough information for our purposes.

A short search for solutions to $\Delta_7(t)=\Delta_9(u)$ modulo squares
(and such that the discriminants do not vanish) leads to the
solution $(t,u)=(-16/3,4)$.
PARI shows that the corresponding cubic fields are both
isomorphic to the unique cubic field of discriminant $-2964$.
(Uniqueness can be seen from the tables
obtainable by ftp at \verb+megrez.math.u-bordeaux.fr+ in
directory \verb+/pub/numberfields+.)
Hence we find a genus-$2$ curve whose Jacobian has a rational
torsion point of order~$63$.
Following the recipe given by Proposition~\ref{explicit}
gives an explicit model for this genus-$2$ curve.
After a few simplifying changes of variable, we obtain the model
\begin{equation}
\label{63torsion}
	C: y^2 = 897x^6 - 197570x^4 + 79136353x^2 - 146398496.
\end{equation}
Let $D$ be the divisor $(R)+(R')-(\infty^+)-(\infty^-)$ on $C$, where 
$$R=\left(\frac{-69+\sqrt{4369}}{2},4515015-68241\sqrt{4369}\right),$$
where $R'$ is the Galois conjugate of $R$, and where $\infty^+$ and $\infty^-$ are
the two points at infinity on a desingularized model of $C$. 
One can check that $D$ maps to a $9$-torsion point on one of the elliptic 
quotients of $C$ and to a $7$-torsion point on the other elliptic quotient
(see~\cite{howeleprevostpoonen}),
so $D$ represents a divisor of order at least $63$ on $C$.
Since $C$ has good reduction at~$5$, and since there is only one
positive multiple of~$63$ less than the Hasse-Weil bound $(1+\sqrt{5})^4$
for $\#(\Jac C)(\F_5)$, we must have $\#(\Jac C)(\F_5)=63$, and hence
the torsion subgroup of $(\Jac C)(\Q)$ is isomorphic to $\Z/63\Z$
and is generated by the class of $D$.
It seems likely that there will be only finitely many genus-$2$ curves over
$\Q$ whose Jacobians possess a rational $63$-torsion point.
It is perhaps even the case that the curve~(\ref{63torsion}) is the only one.

Similarly, we find the solution $(t,u)=(7,-14/13)$ to
$\Delta_7(t)=\Delta_7(u)$ modulo squares.
(We must be careful to exclude solutions where $u=t$, $u=1/(1-t)$,
or $u=(t-1)/t$, since these correspond to taking the same elliptic
curve but with one $7$-torsion point a multiple of the other.)
For these values of $t$ and $u$, the corresponding cubic fields turn out
to be isomorphic, so we indeed obtain a curve $C$  whose Jacobian contains a subgroup
of rational points isomorphic to $\Z/7\Z \times \Z/7\Z$.
Using Proposition~\ref{explicit}, we find the model
$$C: y^2 = x^6 + 3025x^4 + 3232987x^2 + 869675859$$
for this curve.  The Jacobian of the reduction of $C$ modulo $5$
is isogenous to the product of two elliptic curves each with exactly $7$ points
($7$ being the only multiple of $7$ less than the Weil bound), so the
Jacobian of the reduction has exactly $49$ points.
Thus we find that the rational torsion on 
the Jacobian of $C$ is in fact isomorphic to $\Z/7\Z \times \Z/7\Z$.

To handle some of the other cases (in particular those with $N=3$)
we will use the following lemma.
The restrictions on $E$ are not necessary, but we only need the
result under these restrictions.

\begin{lemma}
\label{add3torsion}
Let $k$ be a field whose characteristic is neither $2$ nor $3$.
If $E$ is an elliptic curve over $k$ 
such that $E[2](k)$ is trivial and $j(E) \not=0, 1728$, then there
is a $1$-parameter family of elliptic curves $E'$ over $k$ such that
$E'$ has a $k$-rational $3$-torsion point and $E'[2] \isom E[2]$ as
$\Gal(\kbar/k)$-modules.
\end{lemma}

\begin{proof}
Write $E$ in Weierstrass form as $y^2=x^3+Ax+B$ (so $A,B \not=0$),
and let $r$ be a root of $x^3+Ax+B$.
We claim that specializing $t$ to $-B^2/A^3$ in the universal
elliptic curve $y^2=x^3+(x+t)^2$ over $X_1(3)$ gives one $E'$
with the desired properties.
A calculation shows that $s=-(B/Ar)^2$ is a root of
$x^3+(x-B^2/A^3)^2$ in $k(r)$, and $s$ cannot be in $k$,
since $r$ is at most quadratic over $k(s)$.
Thus $k(r)$ and $k(s)$ are the same cubic extension of $k$,
and hence the curves
	$$y^2=x^3+Ax+B	\qquad \text{and} \qquad y^2=x^3+(x-B^2/A^3)^2$$
have isomorphic $2$-torsion as Galois-modules.

The set of all $E'$ with the desired properties corresponds
to the $k$-rational points of a twist of the modular curve $X_1(2,6)$
classifying elliptic curves with a $3$-torsion point and full level-$2$ structure.
But this modular curve is rational, and the previous paragraph
shows that our twist of it has a $k$-rational point, so our twist
is a rational curve over $k$, and we obtain the desired $1$-parameter family.
\end{proof}

Applying the lemma with $k=\Q(t)$ and $E$ as the universal elliptic
curve with a $7$-torsion point yields a $2$-parameter family of pairs
of elliptic curves with a $3$-torsion point and a $7$-torsion point
respectively, and having isomorphic $2$-torsion as Galois-modules,
so by Proposition~\ref{isotropic}, we obtain a $2$-parameter family of
genus-$2$ curves over $\Q$ whose Jacobians possess a rational $21$-torsion point.
Similarly, if we take $E=E_9^t$ we obtain a $2$-parameter family
of genus-$2$ curves over $\Q$ whose Jacobians have torsion subgroup
containing $\Z/3\Z \times \Z/9\Z$.

Next we construct infinitely many genus-$2$ curves with a
rational $35$-torsion point.
Let $E$ be the elliptic curve $E_7^{-1}$ with a rational $7$-torsion point.
The elliptic curves $E'$ over $\Q$ equipped with a rational
$5$-torsion point and a Galois-module isomorphism $E'[2] \rightarrow E[2]$
correspond to the rational points on a twist $X'$ of $X_1(2,10)$.
Now $X_1(2,10)$ is a covering of $X_1(5)$ with Galois group
$GL_2(\Z/2\Z) \isom S_3$, and the subgroup $A_3$ corresponds
by Galois theory to an intermediate covering whose function field
is the quadratic extension of $\Q(t)$ (where $t$ is the parameter
on $X_1(5)$) obtained by adjoining the square root of the discriminant
of the cubic $f_t(x)$ if $E_5^t$ is written in the form $y^2=f_t(x)$.
This function field is of genus~$1$, since from Table~\ref{discriminants},
$\Delta_5(t)=t(t^2-11t-1)$.
But $X_1(2,10)$ is an elliptic curve as well (see~\cite{kenkumomose}),
so its map down to the intermediate covering must be an isogeny (in fact, a $3$-isogeny).
Similarly our twist $X'$ of $X_1(2,10)$ is a genus-$1$ curve with a
degree-$3$ map to the unique intermediate covering $X''$ of degree~$2$
over $X_1(5)$.
The curve $X''$ classifies elliptic curves $E'$ with a $5$-torsion point
and a Galois-stable $A_3$-orbit of isomorphisms $E'[2] \rightarrow E[2]$.
There are two $A_3$-orbits, and they are defined over
$\Q(\sqrt{\Delta_7(-1)\Delta_5(t)})$, since an automorphism of this
field over $\Q(t)$ is trivial on $\Q(\sqrt{\Delta_7(-1)})$ if and only
if it is trivial on $\Q(\Delta_5(t))$, which means the signatures of
its permutation actions on the nonzero $2$-torsion points of $E$ and
$E_5^t$ must be the same.
Thus $X''$ is the genus-$1$ curve $y^2 = \Delta_7(-1) \Delta_5(t)$, i.e.,
	$$y^2 = -26t(t^2-11t-1).$$
This is an elliptic curve of conductor~$54080$, which is too large for
it to be listed in the tables of~\cite{cremona}, but Cremona's
rank-computing program shows that it has rank~$2$; the points $(-2/13,22/13)$
and $(-26,806)$ are independent of one another and have infinite order.
The genus-$1$ twist $X'$ of $X_1(2,10)$ has a rational point, because
a PARI search finds an elliptic curve $E_5^t$, with $t=1/26$, such
that the cubic field (of discriminant $-104$) obtained by adjoining
a $2$-torsion point is the same as that obtained by adjoining a $2$-torsion
point of $E$.
Thus $X'$ is an elliptic curve 3-isogenous to $X''$.
In particular, $X'$ has rank~$2$, so it has infinitely many rational points,
all but finitely many of which give rise to genus-$2$ curves whose
Jacobians possess a subgroup $\Z/5\Z \times \Z/7\Z \isom \Z/35\Z$.

Similarly the elliptic curve $E=E_9^{-5}$ with a rational $9$-torsion
point has $2$-torsion subgroup isomorphic as Galois-module to that of
the elliptic curve $E_5^{93/10}$, and the elliptic curve
	$$y^2= \Delta_9(-5) t(t^2-11t-1)$$
of conductor $13838400$ (!) has rank~$2$ again according to Cremona's program,
with $(-10/93,6970/93)$ and $(-640/27,5860240/81)$ as
independent points of infinite order.
Thus we obtain infinitely many genus-$2$ curves whose Jacobians possess a
rational $45$-torsion point.

\subsection{Gaining $2$-power torsion}
\label{gaining}

Let $k$ be a field of characteristic not~$2$, let $K$ be a
separable closure of $k$, and let $G_k=\Gal(K/k)$.
If $E$ is an elliptic curve over $k$, then $E[2]\setminus\{0\}\cong\Spec L$
where $L$ is a separable $k$-algebra of dimension~$3$.
Explicitly, if $E$ is in the form $y^2=f(x)$ with $f(x) \in k[x]$
a cubic polynomial, then we may take $L=k[T]/(f(T))$.
As is well known (see~\cite{brumerkramer}, \cite{casselsgenus2}, \cite{schaefer}),
	$$H^1(G_k,E[2]) \isom \ker\left(L^\ast/L^{\ast 2} \stackrel{\Norm}\longrightarrow k^\ast/k^{\ast 2} \right)$$
and the composition
	$$E(k)/2E(k) \hookrightarrow H^1(G_k,E[2]) \isom \ker\left(L^\ast/L^{\ast 2} \stackrel{\Norm}\longrightarrow k^\ast/k^{\ast 2} \right)$$
is a map $\iota$ sending a rational non-$2$-torsion point $P$ with
$x$-coordinate $x_P$ to the image of $x_P-T$.
When $P$ is a non-trivial rational $2$-torsion point, $x_P-T$
vanishes in exactly one component of $L$, and $\iota(P)$ equals
the image of $x_P-T$ is all but this component; the image of $P$
in this last component (which is a copy of $k$) can be deduced up to
squares from knowing that $\iota(P)$ is in the kernel of the norm.

\begin{prop}
\label{gaining2torsion}
Let $f(x)$ and $g(x)$ be cubic polynomials in $k[x]$ such that
	$$E: y^2=f(x) \qquad \text{and} \qquad F: y^2=g(x)$$
are elliptic curves admitting an isomorphism of $G_k$-modules
$\psi: E[2](K) \rightarrow F[2](K)$.
Define $L$ and $\iota$ as above for $E$, and similarly define
$L'$ and $\iota'$ for $F$.
The map $\psi$ induces an isomorphism $\tilde{\psi}: L' \rightarrow L$.
Let $A$ be the quotient of $E \times F$ by the graph of $\psi$.

{\rm (a)} If a point $(P_0,Q_0)$ of $(E \times F)(K)$ maps
to a $k$-rational point on $A$, then $2P_0 \in E(k)$ and $2Q_0 \in F(k)$.

{\rm (b)} Given $P \in E(k)$ and $Q \in F(k)$, there exists a
point $(P_0,Q_0)$ of $(E \times F)(K)$ that maps to
a $k$-rational point on $A$ and such that $2P_0=P$ and $2Q_0=Q$,
if and only if $\iota'(Q)$ corresponds to $\iota(P)$ {\rm(}up to squares{\rm)}
under the isomorphism $\tilde{\psi}$.
\end{prop}

\begin{proof}
Let $\lambda$ be the principal polarization of $A$ derived from the
principal polarization on $E \times F$.
If we compose the isogeny $E \times F \rightarrow A$ with
$\lambda$ and the dual isogeny
$\widehat{A} \rightarrow E \times F$, the result is multiplication-by-2 on
$E \times F$, so part (a) is clear.

Now let $P \in E(k)$ and $Q \in F(k)$.
Suppose that there exists $(P_0,Q_0) \in (E \times F)(K)$ that
maps to a $k$-rational point on $A$ and such that $2P_0=P$ and $2Q_0=Q$.
This means that $(P_0,Q_0)^\sigma - (P_0,Q_0)$ is in the graph of
$\psi$ for all $\sigma \in G_k$.
In particular, under the map induced by $\psi$, the class of
$\xi_\sigma:=P_0^\sigma-P_0$ in $H^1(G_k,E[2])$ is mapped to
the class of $\xi_\sigma':=Q_0^\sigma-Q_0$.
In other words, $\tilde{\psi}$ takes $\iota'(Q)$ to $\iota(P)$.

Conversely suppose that $\tilde{\psi}$ takes $\iota'(Q)$ to $\iota(P)$.
This means that the map induced by $\psi$ takes the image of $P$
in $H^1(G_k,E[2])$ (under the coboundary map) to the image of
$Q$ in $H^1(G_k,F[2])$.
Fix $P_1 \in E(K)$ such that $2P_1=P$ and $Q_1 \in F(K)$
such that $2Q_1=Q$.
Then there exist $2$-torsion points $R \in E[2](K)$ and
$S \in F[2](K)$ such that
	$$\psi(P_1^\sigma-P_1 + R^\sigma-R) = Q_1^\sigma-Q_1 + S^\sigma-S$$
for all $\sigma \in G_k$.
Let $P_0=P_1+R$ and $Q_0=Q_1+R$.
Then $2P_0=P$, $2Q_0=Q$, and
	$$(P_0,Q_0)^\sigma-(P_0,Q_0)$$
is in the graph of $\psi$ for all $\sigma$, so $(P_0,Q_0)$ maps to a
rational point on $A$.
\end{proof}

For elliptic curves $E$ over $\Q$ with all $2$-torsion points rational,
$L$ is simply $\Q \times \Q \times \Q$, the factors corresponding to
the non-trivial torsion points $T_1$, $T_2$, $T_3$.
Now assume that $E=E_{2,8}^t$.
Then $\iota(T_2) \in L$ is $(x_{T_2}-x_{T_1},\ast,x_{T_2}-x_{T_3})$,
where $\ast$ is determined by the condition that the product of all
three components equal~$1$ (modulo squares).
By the formulas in Table~\ref{explicit2torsion}, we have (modulo squares)
	$$\iota(T_2) = \left( -1,-(8t^2-1)(8t^2+8t+1),(8t^2-1)(8t^2+8t+1) \right) \in \left(\Q^\ast/\Q^{\ast 2}\right)^3.$$
Now if $F=E_{2,8}^u$ with non-trivial $2$-torsion points $T_1'$, $T_2'$,
$T_3'$ and corresponding map $\iota'$, then
	$$\iota'(T_2') = \left( -1,-(8u^2-1)(8u^2+8u+1),(8u^2-1)(8u^2+8u+1) \right) \in \left(\Q^\ast/\Q^{\ast 2}\right)^3.$$
If $\psi$ is the isomorphism $E[2] \rightarrow F[2]$ taking $T_i$
to $T_i'$, then the map $\tilde{\psi}$ of Proposition~\ref{gaining2torsion}
is simply the identity
$\Q \times \Q \times \Q \rightarrow \Q \times \Q \times \Q$.
Thus by Proposition~\ref{gaining2torsion}, if there exists $y$ such that
\begin{equation}
\label{128surface}
	(8t^2-1)(8t^2+8t+1) y^2 = (8u^2-1)(8u^2+8u+1),
\end{equation}
then there exists $(P_0,Q_0)$ on $E \times F$ with double $(T_2,T_2')$
such that $(P_0,Q_0)$ maps to a rational point on the quotient $A$
of $E \times F$ by the graph of $\psi$.
In this case, $(P_0,Q_0)$ maps to a new rational $2$-torsion point on $A$,
not in the image of $E(\Q) \times F(\Q)$.

We can consider~(\ref{128surface}) as a genus-$1$ curve over $\Q(t)$,
and we make it an elliptic curve by choosing $(u,y)=(t,1)$ as the origin.
Then the point $(u,y)=(t,-1)$ has infinite order, since the
divisor $(t,-1)-(t,1)$ corresponds to a non-constant point on the Jacobian of
	$$y^2=(8u^2-1)(8u^2+8u+1),$$
which is isomorphic to~(\ref{128surface}) over
$\Q(t)(\sqrt{(8t^2-1)(8t^2+8t+1)})$.
Hence~(\ref{128surface}) is a positive rank elliptic surface
whose points (outside a $1$-dimensional set)
parameterize a family of genus-$2$ curves whose Jacobians
have torsion subgroup over $\Q$ containing
$\Z/2\Z \times \Z/2\Z \times \Z/4\Z \times \Z/8\Z$.

Let us now try to do the same for $E=E_{2,6}^t$ and $F=E_{2,6}^u$.
In this case, from Table~\ref{explicit2torsion} we compute
	$$\iota(T_1) = \left(2(t-3)(t+3)(t-5),(t+3)(t-5),2(t-3)\right) \in \left(\Q^\ast/\Q^{\ast 2}\right)^3.$$
This time in order to get an extra $2$-torsion point on $A$ coming from
a half of $(T_1,T_1')$, we need to find rational solutions to the system
\begin{align}
\label{new72}
		2(t-3)		&= 2(u-3)y^2		\\
\nonumber	(t+3)(t-5)	&= (u+3)(u-5)z^2.
\end{align}
(Note that the third condition
	$$2(t-3)(t+3)(t-5) = 2(u-3)(u+3)(u-5) \;\text{(modulo squares)}$$
would then be automatic.)
If we solve for $t$ in the first equation and substitute into the second,
we obtain the equation
	$$\left( (u-3)y^2+6 \right) \left( (u-3)y^2-2 \right) = (u+3)(u-5)z^2,$$
which defines a genus-$1$ curve over $\Q(u)$.
We make it an elliptic curve by choosing $(y,z)=(1,1)$ as origin,
and then note that $(y,z)=(-1,1)$ is a point of infinite order, since it is of infinite order for the specialization $u=0$.
Thus the system~(\ref{new72}) provides us with a positive rank
elliptic surface whose points parameterize a family of
genus-$2$ curves whose Jacobians have torsion subgroup over $\Q$
containing $\Z/2\Z \times \Z/6\Z \times \Z/6\Z$.

Next we investigate the possibility of gaining $2$-power torsion
when $E=E_{2,6}^t$ and $F=E_{2,8}^u$.
Let $T_1,T_2,T_3$ and $T_1',T_2',T_3'$ be the nontrivial $2$-torsion
points on $E$ and $F$, respectively, as in Table~\ref{explicit2torsion}.
We have
\begin{align*}
	\iota(T_1)	& = \left(2(t-3)(t+3)(t-5),(t+3)(t-5),2(t-3)\right) \in \left(\Q^\ast/\Q^{\ast 2}\right)^3,	\\
	\iota'(T_2')	& = \left( -1,-(8u^2-1)(8u^2+8u+1),(8u^2-1)(8u^2+8u+1) \right) \in \left(\Q^\ast/\Q^{\ast 2}\right)^3.
\end{align*}
In an attempt to obtain simpler equations than we would by
mapping $T_i$ to $T_i'$ for each $i$, we let $\psi:E[2] \rightarrow F[2]$
be the isomorphism such that $\psi(T_1)=T_2'$, $\psi(T_2)=T_3'$,
and $\psi(T_3)=T_1'$.
Hence $\tilde{\psi}$ is the isomorphism
$\Q \times \Q \times \Q \rightarrow \Q \times \Q \times \Q$
acting on the factors as the permutation $(1\;3\;2)$.
By Proposition~\ref{gaining2torsion}, a point on $E \times F$
with double $(T_1,T_2')$ maps to a new rational point on the
quotient $A$ if and only if we can find rational numbers $y$
and $z$ such that
\begin{align*}
\label{new96}
	2(t-3)	&= (-1)y^2,	\\
\nonumber	(t+3)(t-5)	&= (8u^2-1)(8u^2+8u+1) z^2.
\end{align*}
If we solve the first equation for $t$, substitute into the second,
and multiply both sides by~$4$, we obtain
	$$(y^2-12)(y^2+4) = 4 (8u^2-1)(8u^2+8u+1) z^2.$$
For $y=2/9$, the resulting genus-$1$ curve has a rational point
$(u,z)=(1/3,44/9)$, and is birational to the elliptic curve
	$$Y^2=X^3-1681X$$
of conductor $53792$ and rank~$2$.
The new rational point on $A$ is a $2$-torsion point, since its double
is the image of $(T_1,T_2')$, which is in the graph of $\psi$.
Hence we have produced a family of genus-$2$ curves,
parameterized by the points on a positive rank elliptic curve,
whose Jacobians have torsion subgroup over $\Q$ containing
$\Z/2\Z \times \Z/2\Z \times \Z/24\Z$.

\section{Genus three}
\label{genus3sect}
\subsection{Jacobians $(2,2,2)$-isogenous to a product of elliptic curves}
\label{genus3}

In this section we will show how one can find a curve of genus~$3$ whose Jacobian
is isogenous over a quadratic extension of the base field
to a product of three given elliptic curves.  
Genus-$3$ curves of the sort we will see were used in \cite{casselsquartic}.

We maintain the conventions of Section~\ref{conventions}.

Suppose $E_1$, $E_2$, and $E_3$ are elliptic curves over a separably closed 
field $K$, and let $N$ be a positive integer not divisible by the
characteristic of $K$.
The product of the canonical polarizations on the $E_i$ is a 
principal polarization $\lambda$  on the abelian variety $A=E_1\times E_2\times E_3$,
and the Weil pairings on the $N$-torsion subgroups of the $E_i$ combine
to give us a non-degenerate alternating pairing $e_N$ from $A[N]$ to the
group scheme of $N$th roots of unity over $K$.
Suppose $G$ is a sub-group-scheme of $A[N]$ that is maximal isotropic
with respect to the pairing $e_N$.  As in the similar situation we saw in
Section~\ref{construction2}, the polarization $N\lambda$ on $A$ reduces to
a principal polarization
$\mu$ on the quotient variety $B=A/G$.  A result of Oort and Ueno
\cite{oortueno} shows that the polarized variety $(B,\mu)$ either
breaks up as a product of lower-dimensional polarized varieties or is 
the canonically polarized Jacobian
of a curve $C$ over $K$ of genus~$3$.\footnote{
	The analogous statement is not necessarily true over a field
	that is not separably closed.  See the remark following the
	proof of Proposition~\ref{quartic3}.}
We would like to see which group-schemes $G$ lead to curves in the case
where $N=2$.  Since we will be working over a separably closed field,
we will identify sub-group-schemes of $A[2]$ with subgroups of $A[2](K)$.

\begin{lemma}
\label{genus3isotropic}
Let $A=E_1\times E_2\times E_3$ and $e_2$ be as above.  
There are exactly $135$ maximal isotropic subgroups $G$ of $A[2](K)$.
Exactly $81$ of these group-schemes are of the form $G_1\times G_2$,
where $G_1$ is a maximal isotropic subgroup of $E_i[2](K)$ for some $i$
and $G_2$ is a maximal isotropic subgroup of $\prod_{j\neq i} E_j[2](K)${\rm;}
for these $G$, the polarized variety $(B,\mu)$ splits into a product of
lower-dimensional polarized varieties.
If $G$ is one of the remaining $54$ groups,
then for each $i$ we may label the nonzero elements of $E_i[2](K)$ by the
symbols $P_i$, $Q_i$, and $R_i$ in such a way so that $G$ is the group
\begin{align*}
\big\{ (0,0,0), &(0,Q_2,Q_3), (Q_1, 0, Q_3), (Q_1, Q_2, 0),\\
     &(P_1, P_2, P_3), (P_1, R_2, R_3), (R_1, P_2, R_3) , (R_1, R_2, P_3)\big\}.
\end{align*}
\end{lemma}

\begin{rem}
In fact, our constructions below will show that every group of the last type
gives rise to a curve.  This fact can also be proven by assuming that a $G$
of the given type is the kernel of a map $A\rightarrow A_1\times A_2$
of polarized varieties and obtaining a contradiction. In anticipation
of this result, we will call maximal isotropic subgroups of $A[2](K)$
(or sub-group-schemes of $A[2]$) {\em non-split\/} if they are of the
latter type.
\end{rem}

\begin{proof}
In any group isomorphic to $(\Z/2\Z)^6$ with a non-degenerate alternating pairing,
there are $(2^6-1)(2^5-2)(2^4-2^2)$ ways of choosing an ordered triple
$(v_1,v_2,v_3)$ that generate a maximal isotropic subgroup,
and each maximal isotropic subgroup has $(2^3-1)(2^3-2)(2^3-2^2)$
such bases, so there are
	$$\frac{(2^6-1)(2^5-2)(2^4-2^2)}{(2^3-1)(2^3-2)(2^3-2^2)} = 135$$
such subgroups.

For $i=1,2,3$ let $S_i$ denote the set of maximal isotropic subgroups
of $A[2](K)$ that can be written $G_1\times G_2$, with $G_1\subset E_i[2](K)$
and $G_2\subset \prod_{j\neq i} E_j[2](K)$.  There are $3$ choices for $G_1$
and $(2^4-1)(2^3-2)/(2^2-1)(2^2-2)=15$ choices for $G_2$, so $\#S_i = 45$.
The intersection of any two of the $S_i$ is the set $S$ of subgroups
of $A[2](K)$ that can be written $G_1\times G_2\times G_3$, with $G_i\subset E_i[2](K)$;
clearly $\#S = 27$.  Thus there are $\#S_1 + \#S_2 + \#S_3 - 2\#S = 81$
subgroups that split as in the statement of the lemma.

We are left with $135-81 = 54$ subgroups to account for, and it is easy to see
that there are exactly this many subgroups of the form described in the
final sentence of the lemma:  There are $3$ choices for each of the
$Q_i$, and given $P_1$ and $P_2$, there are $2$ choices for $P_3$.
\end{proof}

Suppose now that $k$ is a field of characteristic not~$2$
with separable closure $K$, and let $E_1$,
$E_2$, and $E_3$ be elliptic curves over $k$.  It is clear that a non-split
sub-group-scheme of $E_{1K}\times E_{2K}\times E_{3K}$ will come from a
sub-group-scheme of $E_1\times E_2\times E_3$ if and only if we have 
both that all of the points $Q_i$ are defined over $k$ and
that every $k$-automorphism of $K$ that moves any of the $P_i$ moves
exactly two of them.  Given the first condition, the second condition
will hold if and only if the product of the discriminants of the curves $E_i$
is a square in $k$.  

So suppose $E_i$ is the elliptic curve over $k$
given by the equation $y^2 = x(x^2 + A_i x + B_i)$ where
$B_i\neq 0$, and let $Q_i$ be the rational $2$-torsion point $(0,0)$.
The discriminant of $E_i$ differs by a square factor from
the number $\Delta_i = A_i^2 - 4B_i$,
so let us assume that $\Delta_1\Delta_2\Delta_3$ is a square in $k$.
For each $i$ let $P_i$ be a nonzero element of $E_i[2](K)$ different from $Q_i$.  
Let $G$ be the non-split sub-group-scheme of the product 
$A=E_1\times E_2\times E_3$ corresponding to this choice of $P$'s and $Q$'s.
We are led to the question: Is the quotient polarized variety $A/G$ 
the Jacobian of a curve over $k$, and if so, what equations define the curve?

Before we answer this question, we must define some numbers.
For each $i$, we let $d_i = -(A_i + 2 x_{P_i})$, where $x_{P_i}$ denotes 
the $x$-coordinate of the point $P_i$.  Note that $d_i^2 = \Delta_i$.
The product  $R=d_1 d_2 d_3$ is an element of $k$ because 
$\Delta_1\Delta_2\Delta_3$ was assumed to be a square.
We let $\lambda_i$ denote $A_i/d_i$, and 
we define the {\em twisting factor\/} (associated to the given $E_i$, 
$P_i$, and $Q_i$) to be the number
\begin{align*}
T  & = R\left( \frac{A_1^2}{\Delta_1}
             +\frac{A_2^2}{\Delta_2}
             +\frac{A_3^2}{\Delta_3} 
             - 1 \right)  - 2 A_1 A_2 A_3 \\
   & = d_1 d_2 d_3 \left(\lambda_1^2 + \lambda_2^2 + \lambda_3^2
                        - 2\lambda_1\lambda_2\lambda_3 - 1\right).
\end{align*}
(The twisting factor is so named because it determines a quadratic extension 
of $k$ over which $A/G$ becomes isomorphic to a Jacobian; see Proposition~\ref{quartic3}.)

\begin{prop}
\label{hyperelliptic3}
With notation as above, suppose $T=0$.  Then each of the products
$B_1B_2$, $B_1B_3$, and $B_2B_3$ is a square, and $A/G$ is isomorphic
{\rm(}over $k${\rm)} to the polarized Jacobian of the hyperelliptic curve $C$ over $k$
defined by the homogeneous equations
\begin{align*}
W^2Z^2 &= aX^4 + bY^4 + cZ^4\\
    0  &= dX^2 + eY^2 + fZ^2,
\end{align*}
where $a$, $b$, and $c$ are given by
\begin{align*}
 a &= \left(\frac{RB_1}{2}\right) 
      \left(         -  \frac{B_1}{\Delta_1} 
                     +  \frac{B_2}{\Delta_2} 
                     +  \frac{B_3}{\Delta_3}  \right) \\
 b &= \left(\frac{RB_2}{2}\right) 
      \left(\phantom{-} \frac{B_1}{\Delta_1}
                     -  \frac{B_2}{\Delta_2}
                     +  \frac{B_3}{\Delta_3}  \right) \\
 c &= \left(\frac{RB_3}{2}\right) 
      \left(\phantom{-} \frac{B_1}{\Delta_1}
                     +  \frac{B_2}{\Delta_2}
                     -  \frac{B_3}{\Delta_3}  \right),
\end{align*}
where $d$, $e$, and $f$ are determined up to sign by the relations 
\begin{align*}
 B_2B_3d^2 &= 1 \\
 B_1B_3e^2 &= 1 \\
 B_1B_2f^2 &= 1,
\end{align*}
and where the signs of $d$, $e$, and $f$ are chosen so that we have
$A_1 = - aef$ and $A_2 = - bdf$ and $A_3 = - cde$.
\end{prop}

\begin{proof}
The statement that $T=0$ is equivalent to the statement that 
$\lambda_1^2+\lambda_2^2+\lambda_3^2  - 2 \lambda_1\lambda_2\lambda_3 - 1 = 0.$
Solving this equation for $\lambda_3$ in terms of $\lambda_1$
and $\lambda_2$ leads to
$$\lambda_3 = \lambda_1\lambda_2 \pm \sqrt{(\lambda_1^2 - 1)(\lambda_2^2 - 1)},$$
and dividing this last equality by $d_3$ gives
$$\frac{A_3}{\Delta_3} = \frac{A_1 A_2 \pm 4\sqrt{B_1 B_2}}{R}.$$
Thus  $B_1B_2$ is a square in $k$.  By symmetry, the
numbers $B_1B_3$ and $B_2B_3$ are squares as well.

We leave it to the reader to verify that the signs of $d$, $e$, and $f$ can 
be chosen so that the relations $A_1 = - aef$ and $A_2 = - bdf$ and $A_3 = - cde$
hold; this can be seen by noting that the squares of the
desired relations, as well as the product of the desired relations,
follow from the formulas given and the condition that $T=0$.

Note that the coefficients $d$, $e$, and $f$ are all nonzero, so 
the single equation $dX^2+eY^2+fZ^2 = 0$ defines a nonsingular curve $C'$
of genus~$0$.
The fact that none of the $B_i$ is zero implies that
none of the numbers $a/d^2 + b/e^2$, $a/d^2 + c/f^2$, and $b/e^2 + c/f^2$ is zero.
Furthermore, the nonvanishing of the $B_i$ together with 
the fact that $T$ is zero can be used to show that
$$\Big(\frac{a}{d^2}\Big)\Big(\frac{b}{e^2}\Big)
+ \Big(\frac{a}{d^2}\Big)\Big(\frac{c}{f^2}\Big)
+ \Big(\frac{b}{e^2}\Big)\Big(\frac{c}{f^2}\Big)$$
is not zero.  We leave it to the reader to show that the nonvanishing of
the several expressions just mentioned guarantees that the
element $a(X/Z)^4 + b(Y/Z)^4 + c$ of the function field of $C'$
has exactly $8$ zeroes, each simple, and exactly $2$ poles, each of order $4$.
Since the total quotient ring of $C$ is obtained from the function field of $C'$
by adjoining a square root of this function, we see 
that $C$ is a hyperelliptic curve of genus~$3$.

It will suffice to prove that $\Jac C \cong A/G$ in the special case where 
$k$ has characteristic $0$, for if $k$ has positive characteristic we can simply lift
all of the coefficients $A_i$ and $B_i$ up to the ring $W$ of Witt vectors over~$k$;
to see that this can be done in such a way that $T$ lifts to $0$, we argue as
follows. 
First we lift each $\Delta_i$ up to $W$ in such a way that the product of
the lifted values is a square in $W$, and we lift $R$ to a square root of
this product.  Now we view $T$ as a function of the three variables $A_i$.
We will be able to use Hensel's lemma to lift the $A_i$ up to $W$ so as to
make $T=0$ if any one of the partial derivatives $\partial T/\partial{A_i}$
is nonzero.  We claim that at least one of these derivatives is nonzero.
To prove this, let us assume that all three of the partial derivatives are
zero and obtain a contradiction.  From our assumption we find that 
$RA_i = \Delta_i\prod_{j\neq i}A_j$ for each $i$. 
If any one of the $A_i$ were zero, these three equalities would imply that
all of the $A_i$ were zero, which would contradict the assumption that $T=0$.  
But if all of the $A_i$ were nonzero,
then by multiplying the three equalities together we would find that $R=A_1A_2A_3$,
and this formula for $R$, combined with 
$RA_i = \Delta_i\prod_{j\neq i}A_j$, would show that $A_i^2 = \Delta_i$, which 
would lead to the impossibility $B_i = 0$.  This proves our claim.

To prove that $\Jac C \cong A/G$ in characteristic zero
we need only consider the universal case, in which we let 
the $d_i$ and the $\lambda_i$ be indeterminates, we let 
$\ell$ be the quotient field of the domain
$$\Q[d_1, d_2, d_3, \lambda_1,\lambda_2,\lambda_3]/
   (\lambda_1^2 + \lambda_2^2 + \lambda_3^2 - 2\lambda_1\lambda_2\lambda_3-1),$$
we take
\begin{align*}
A_i &= \lambda_i d_i \\
B_i &= (A_i^2 - d_i^2)/4 \\
Q_i &= (0,0) \\
P_i &= (-(d_i+A_i)/2, 0) \\
  R &= d_1 d_2 d_3,
\end{align*}
and we let $k$ be the subfield $\Q(A_1, A_2, A_3, B_1, B_2, B_3, R)$ of $\ell$.
(Note that $k$ is fixed by the involution of $\ell$ that acts on
$\lambda_1$, $\lambda_2$, $d_1$, and $d_2$ by multiplication by $-1$,
so $k$ is a proper subfield of $\ell$ and contains (by symmetry) none
of the $d_i$.)  So let us assume that we are in the universal case.

Let $\iota_X$ be the involution on $C$ that fixes $W$, $Y$, and $Z$ and that
sends $X$ to~$-X$.
The involution $\iota_X$ gives us a double cover 
$\varphi_X\colon C\rightarrow F_X$ of curves over $k$,
and we would like to find equations for the curve $F_X$.
If we dehomogenize the equations for $C$ with respect to $Z$
by letting $w=W/Z$, $x=X/Z$, and $y=Y/Z$,
and if we then divide by $\iota_X$ by defining $u=x^2$, we find that
the quotient curve $F_X$ is given by
\begin{align*}
w^2 &= au^2 + by^4 + c\\
 0  &= du + ey^2 + f.
\end{align*}
This pair of equations can be combined to get the single equation
\begin{equation}
\label{FXequation}
v^2 = (ae^2 + bd^2) y^4 + 2aef y^2 + (af^2 + cd^2),
\end{equation}
where $v=dw$.
Using Example~3.7 (pp.~293--294) in Section~X of \cite{silverman},
we see that the Jacobian of the genus-$1$ curve $F_X$ is the elliptic
curve $E_X$ over $k$ defined by
\begin{equation}
\label{EXequation}
y^2 = x(x^2 + A_X x + B_X)
\end{equation}
where $A_X = -aef$ and $4 B_X = (aef)^2 - (ae^2+bd^2)(af^2+cd^2)$.
Clearly we have $A_X = A_1$, and by using a little algebra and
the fact that $T=0$ we can see that $B_X = B_1$. Thus the double
cover $\varphi_X$ gives us a map $\varphi_X^*$ from $E_1$ to the
Jacobian $J$ of $C$.

Similarly, the involution $\iota_Y$ on $C$ that sends $Y$ to $-Y$
and that fixes the other variables
gives us a degree-$2$ cover $\varphi_Y$ from $C$ 
to the curve $F_Y$ over $k$ given by the equation
$$v^2 = (ae^2 + bd^2) x^4 + 2bdf x^2 + (bf^2 + ce^2).$$
The Jacobian of $F_Y$ is isomorphic to $E_2$, so we get 
a map $\varphi_Y^*$ from $E_2$ to $J$.

Lastly, the involution $\iota_Z = \iota_X\circ\iota_Y$ 
us a degree-$2$ cover $\varphi_Z\colon C\rightarrow F_Z$
to a curve $F_Z$ over $k$ whose Jacobian is isomorphic to $E_3$,
so we get a map $\varphi_Z^*$ from $E_3$ to $J$.
(When dividing $C$ by $\iota_Z$, the reader may find it helpful 
to recast the first defining equation of $C$ into the form
$V^2Y^2 = aX^4+bY^4+CZ^4$ by letting $V=WZ/Y$; this will make
it possible to dehomogenize with respect to $Y$ and get equations
similar to the ones obtained when dividing by $\iota_X$.)

Let $I$ denote the subgroup of the automorphism group of
$J$ generated by $\iota_X^*$, $\iota_Y^*$, and $\iota_Z^*$,
and let $\CC$ denote the category of abelian varieties over
the separable closure $K$ of $k$ up to isogeny.
The semisimple group ring $\Q[I]$ acts on the class $[J]$ of $J$ in $\CC$,
and $[J]$ splits into the direct sum of its eigenspaces.  
The class of the image of $\varphi_X^*$ consists of 
the sum of the eigenspaces on which $\iota_X^*$ acts as $1$,
while the class of the image of $\varphi_Y^*$ consists of 
the sum of the eigenspaces on which $\iota_Y^*$ acts as $1$.
However, the eigenspaces on which both $\iota_X^*$ and $\iota_Y^*$
act as $1$ are trivial, because their sum is the class of the
Jacobian of the quotient of $C$ by the group $\langle\iota_X,\iota_Y\rangle$,
and this quotient has genus $0$.  
Thus the classes of the images of $\varphi_X^*$ and
$\varphi_Y^*$ have trivial intersection.  Similarly, we find that 
the class of the image of $\varphi_Z^*$ shares no nonzero
eigenspaces with the classes of $\varphi_X^*$ or $\varphi_Y^*$.
It follows that the  morphism
$\Phi^* = \varphi_X^*\times\varphi_Y^*\times\varphi_Z^*$ from
$A$ to $J$ is an isogeny. 

Let $\mu$ denote the canonical polarization of $J$.
The fact that $\varphi_X$ has degree $2$ implies that
$\widehat{\varphi_X^*}\mu\varphi_X^*$
is the multiplication-by-$2$ map on $E_1$, and similarly 
$\widehat{\varphi_Y^*}\mu\varphi_Y^*$ and $\widehat{\varphi_Z^*}\mu\varphi_Z^*$
are the multiplication-by-$2$ maps on $E_2$ and $E_3$; 
here $\widehat{\ }$ indicates the dual morphism.
We find that we have a commutative diagram
$$
\begin{matrix}
A                 & \mapright{2\lambda} & \widehat{A} \\
\Bmapdown{\Phi^*} &                     &\Bmapup{\widehat{\Phi^*}}\\
J                 & \mapright{\mu}      & \widehat{J}.
\end{matrix}
$$
We see that $\Phi^*$ has degree eight
and its kernel $G'$ is a maximal isotropic sub-group-scheme of $A[2]$.
The group-scheme $G'$ must be non-split. 
To complete the proof of the proposition 
we must show that it is equal to the given group-scheme $G$.

Since $G'$ is non-split, there must be nonzero elements $Q_i'\in E_i[2](K)$
such that  $G'(K)$ contains the three elements 
$(0, Q_2', Q_3')$, $(Q_1', 0, Q_3')$, and $(Q_1', Q_2', 0)$.
As we noted before the statement of the proposition, the points $Q_i$ must 
actually be defined over $k$.  However, because $d_i$ is not in $k$,
the curve $E_i$ has only one nonzero $k$-defined $2$-torsion point, namely $Q_i$.
Therefore $Q_i' = Q_i$.
There are exactly two non-split maximal isotropic subgroups of $A[2](K)$
that contain the subgroup
$$H=\{(0,0,0), (0, Q_2, Q_3), (Q_1, 0, Q_3), (Q_1, Q_2, 0)\},$$
and the intersection of these two groups is $H$.  So to show that $G'(K)=G(K)$,
all we must do is show that the two groups 
contain a common element that is not in $H$.  To show this, we can specialize
our universal example to a particular case, and show that the specialized
groups $G$ and $G'$ contain a common element not in $H$.

Consider the specialization map $\ell\ra\C$ that takes each $d_i$ to $-4$
and each $\lambda_i$ to $-1/2$; we will abuse notation by saying
that $d_1 = -4$, and so on.  We see that each $A_i = 2$, each $B_i = -3$,
each $\Delta_i = 16$, each $P_i = (1,0)$, 
and the curve $C$ is defined by the two equations
$W^2Z^2 = -18(X^4 + Y^4 + Z^4)$ and $0=(1/3)(X^2 + Y^2 + Z^2)$.

Let $\Phi_*$ denote the map
$\varphi_{X*} \times\varphi_{Y*} \times\varphi_{Z*}$ from
$J$ to $A$ and recall that $\Phi^*$ denotes the map
$\varphi_X^* \times\varphi_Y^* \times\varphi_Z^*$ from $A$ to $J$.
An easy computation shows that $\Phi^*\Phi_*$ is multiplication by $2$ on $J$.
Thus the image under $\Phi_*$ of a $2$-torsion element of $J$ is in the
kernel of $\Phi^*$.  To complete the proof, we will show that 
$(P_1,P_2,P_3)\in A(K)$ is the image under $\Phi_*$ of a $2$-torsion point of $J$.

Let $U_1$ and $U_2$ be the Weierstrass points on $C$ given in homogeneous
coordinates $[W:X:Y:Z]$ by $[0:\zeta:\zeta^2:1]$ and $[0:\zeta^2:\zeta:1]$,
respectively, where $\zeta = e^{2\pi i/3}$.
Note that $U_1 - U_2$ represents a $2$-torsion
element of $J$.  Under the specialization we have made,
equation~(\ref{FXequation}), which defines $F_X$, becomes 
$v^2 = -4 y^4 -4 y^2 -4,$ and in these $(y,v)$ coordinates
we have $\varphi_X(U_1) = (\zeta^2,0)$ and $\varphi_X(U_2) = (\zeta,0)$.
The curve $E_X$, defined by equation~(\ref{EXequation}), is given by
$y^2 = x^3 + 2 x^2 -3x,$ and under the isomorphism
$$(y,v)\mapsto 
  \left(\frac{iv - y^2 - 2}{y^2}, \frac{-2v-2i y^2 - 4i}{y^3}\right)$$ from
$F_X$ to $E_X$ these two points on $F_X$ map to 
$(i\sqrt{3}, 2\zeta\sqrt{3})$ and $(-i\sqrt{3}, 2\zeta^2\sqrt{3})$,
respectively.  The difference of these two points on the elliptic curve
$E_X$ is equal to the $2$-torsion point $(1,0)$.
Thus $\varphi_{X*}(U_1-U_2) = P_1$. By symmetry, we find that
$\varphi_{Y*}(U_1-U_2) = P_2$ and $\varphi_{Z*}(U_1-U_2) = P_3$
as well, so $(P_1,P_2,P_3)$ is in the image of $J[2](K)$ under $\Phi_*$ and hence
in $G'(K)$.  It is in $G(K)$ as well, so $G'=G$ and the proposition is proved.
\end{proof}

\begin{prop}
\label{quartic3}
Let notation be as before Proposition~\ref{hyperelliptic3}, and
suppose $T\neq 0$.  Let $C$ be the plane quartic over $k$ defined by
$$B_1X^4 + B_2Y^4 + B_3Z^4 + dX^2Y^2 + eX^2Z^2 + fY^2Z^2 = 0,$$
where 
\begin{align*}
 d &= \frac{1}{2} \left(-A_1A_2 + \frac{A_3R}{\Delta_3}\right)\\
 e &= \frac{1}{2} \left(-A_1A_3 + \frac{A_2R}{\Delta_2}\right)\\
 f &= \frac{1}{2} \left(-A_2A_3 + \frac{A_1R}{\Delta_1}\right).\\
\end{align*}
Let $k'$ be the field $k(\sqrt{T})$. 
Then the polarized Jacobian of $C_{k'}$ is isomorphic to
the polarized variety $A_{k'}/G_{k'}$.
\end{prop}

\begin{proof}
We leave it to the reader to show that an equation of the form
$$aX^4 + bY^4 + cZ^4 + dX^2Y^2 + eX^2Z^2 + fY^2Z^2 = 0$$
defines a non-singular curve if and only if the seven numbers
$a$, $b$, $c$, $d^2 - 4ab$, $e^2-4ac$, $f^2-4bc$, and 
$af^2 + be^2 + cd^2 - 4abc - def$ are nonzero.
In our case, these numbers are
$B_1$, $B_2$, $B_3$, $TR/4\Delta_3$, $TR/4\Delta_2$, $TR/4\Delta_1$,
and $T^2/16$, which are all nonzero.  Thus our $C$ is a non-singular
curve of genus~$3$.

As in the proof of Proposition~\ref{hyperelliptic3}, we quickly reduce
to the universal case.  This time, that means that the $d_i$ and the
$\lambda_i$ are indeterminates, that $\ell$ is the field 
$$\Q(d_1, d_2, d_3, \lambda_1,\lambda_2,\lambda_3),$$
that
\begin{align*}
A_i &= \lambda_i d_i \\
B_i &= (A_i^2 - d_i^2)/4 \\
Q_i &= (0,0) \\
P_i &= (-(d_i+A_i)/2, 0)\\
  R &= d_1 d_2 d_3,
\end{align*}
and that $k$ is the proper subfield $\Q(A_1, A_2, A_3, B_1, B_2, B_3, R)$ of $\ell$.
Let $\ell'=\ell(\sqrt{T})$;
note that $k'=k(\sqrt{T})$ is a proper subfield of $\ell'$
because it contains none of the $d_i$.

Let $\iota_X$ be the involution $X\mapsto -X$ of $C_{k'}$ and
let $\varphi_X\colon C_{k'}\ra F_X$ be the double cover induced by $\iota_X$.
To find a model for the curve $F_X$ over $k'$,
we dehomogenize the equation for $C$ by letting  $x=X/Z$ and $y=Y/Z$;
then, setting $u=x^2$, we find the model
$$B_1 u^2 + B_2 y^4 + B_3 + duy^2 + eu + fy^2 = 0$$
for $F_X$.
If we let $v=2B_1 u + d y^2 + e$ and simplify, we get the model
$$v^2 = (d^2-4B_1B_2) y^4 + (2de-4B_1f) y^2 + (e^2-4B_1B_3).$$
Example 3.7 (pp.~293--294) of \cite{silverman} shows that
the Jacobian $E_X$ of $F_X$ is the elliptic curve over $k'$
defined by $y^2 = x^3 + A_X x^2 + B_X x$,
where $A_X = 2B_1 f - de$ and 
$B_X = B_1 ( B_1 f^2 + B_2 e^2 + B_3 d^2 - def - 4B_1 B_2 B_3)$.
Using the formulas for $d$, $e$, and $f$ given in the
proposition, we find that $A_X = A_1 T/4$ and $B_X = B_1 T^2/16$.
Thus we see that $E_X \cong E_{1k'}$,
and the double cover $\varphi_X\colon C_{k'}\ra F_X$
gives us a map $\varphi_X^*$ from $E_{1k'}$ to the Jacobian $J$
of $C_{k'}$.

If we define two more involutions $\iota_Y$ and $\iota_Z$ of $C_{k'}$ in the
obvious way, we get double covers $\varphi_Y\colon C_{k'}\ra F_Y$ and
$\varphi_Z\colon C_{k'}\ra F_Z$ that give rise to
homomorphisms $\varphi_Y^*\colon E_{2k'}\ra J$ and
$\varphi_3^*\colon E_{3k'}\ra J$.

Let $\mu$ be the canonical polarization of 
$J$ and let $\lambda$  be the product polarization on 
$A=E_1\times E_2\times E_3$.
As in the proof of Proposition~\ref{hyperelliptic3}, 
we get a commutative diagram
$$
\begin{matrix}
A_{k'}                 & \mapright{2\lambda} & \widehat{A_{k'}} \\
\Bmapdown{\Phi^*} &                     &\Bmapup{\widehat{\Phi^*}}\\
J                 & \mapright{\mu}      & \widehat{J}
\end{matrix}
$$
where $\Phi^* = \varphi_X^*\times\varphi_Y^*\times\varphi_Z^*$
is an isogeny of degree $8$ whose kernel $G'$ is a non-split maximal isotropic
sub-group-scheme of $A_{k'}[2]$.  Our task is to show that $G'=G_{k'}$.

Let $K$ be an algebraic closure of $\ell'$.
Rationality arguments as in the proof of Proposition~\ref{hyperelliptic3}
show that $G'(K)$ contains the subgroup 
$$H=\{(0,0,0), (0, Q_2, Q_3), (Q_1, 0, Q_3), (Q_1, Q_2, 0)\},$$
and, as before, to show that $G'(K)=G(K)$ all we must do is
show that the two groups contain a common element that is not in $H$.
To show this, we once again specialize our universal example to a
particular example.

Consider the specialization map $\ell\ra\C$ that takes each $d_i$ to $4$
and each $\lambda_i$ to $-1/2$; we will abuse notation by saying
that $d_1 = 4$, and so on.  We see that each $A_i = 2$, each $B_i = -3$,
each $\Delta_i = 16$, each $P_i = (-3,0)$, and $T=-32$.  For each $i$,
let $R_i$ be the $2$-torsion point $(1,0)$ of $E_i$.
Note that the $E_i$ are the same as in the specialization at the end
of the proof of Proposition~\ref{hyperelliptic3}.  From that proof,
we know that the polarized quotient
of $A$ by the subgroup generated by $H$ and $(R_1, R_2, R_3)$ is
the Jacobian of a hyperelliptic curve.
Since $A/G'$ is the Jacobian of a plane quartic, Torelli's theorem
shows that $G'(K)$ cannot possibly contain $(R_1,R_2,R_3)$.  The only
possibility remaining is that $G'(K)$ contains $(P_1,P_2,P_3)$,
which shows that $G'=G$ and completes the proof.
\end{proof}

\begin{rem}
One might ask whether the base field extension to $k(\sqrt{T})$
is necessary for the proposition to be true.  Indeed it is necessary.
To see this, consider an arbitrary plane quartic $C$ over $k$,
let $J$ be its polarized Jacobian, and let $K$ be  a separable closure
of $k$.  Since $C$ is not hyperelliptic,
we have an isomorphism $\Aut J \cong \Aut C \times \{\pm1\}$ of
Galois modules (where the Galois action on $\{\pm1\}$ is trivial).  
Taking Galois cohomology, we find 
$$H^1(\Aut J)\cong H^1(\Aut C) \times \Hom(\Gal(K/k),\{\pm1\}).$$
The two $H^1$'s catalog the twists of $J$ and $C$, respectively, and
the $\Hom$ catalogs field extensions of $k$ of degree at most $2$.
Suppose $J'$ is a quadratic twist of $J$ corresponding to
an element of $H^1(\Aut J)$ that is trivial in $H^1(\Aut C)$ but nontrivial
in the $\Hom$.
The curve over $k$ that one obtains from $J'$ is none other than $C$,
and it takes a quadratic extension to make $\Jac C$ isomorphic to $J'$.
\end{rem}

We can use Proposition~\ref{quartic3} to give an example of a Jacobian of
a curve over $\Q$ whose conductor, while not exactly {\em small},
is at least not so big.  Recall that Mestre's result \cite{mestreconducteurs}
implies that under standard conjectures the conductor of a $3$-dimensional
abelian variety over $\Q$ is at least $1100$.

\begin{cor}
\label{conductor2940}
The conductor of the Jacobian of the curve 
$$2X^4 + 2Y^4 + 15Z^4 + 3X^2Y^2 - 11X^2Z^2 -11Y^2Z^2 = 0$$
is $2940$.
\end{cor}

\begin{proof}
Take $E_1$ and $E_2$ to be the curve $y^2 = x^3 -11x^2+32x$, which is isomorphic
to the curve 14A4 of \cite{cremona} and has conductor $14$.
Take $E_3$ to be the curve $y^2 = x^3 - 31x^2+240x$, which is isomorphic to
the curve 15A3 of \cite{cremona} and has conductor $15$.  If we take
$P_1=P_2$ to be a nonzero $2$-torsion point on $E_1$ other than $(0,0)$,
and if we take $P_3$ to be $(15,0)$, then we find that the twisting factor
is $T=32^2$.  Applying Proposition~\ref{quartic3} to these curves gives
the curve in the statement of the corollary.
\end{proof}

\subsection{Building hyperelliptic Jacobians --- introduction}

In the next few sections we will find triples $(E_1, E_2, E_3)$ of elliptic
curves over $\Q$ that have large rational torsion subgroups and
for which we can choose $2$-torsion points $P_i$ and $Q_i$
that make the twisting factor equal to zero.  Our strategy will be
to specify the rational torsion structure on $E_3$ and determine the
corresponding conditions on $E_1$ and $E_2$ that will make the twisting
factor zero. We will not exhaust the possible combinations
of torsion structures; the
equations that arise become very messy very quickly, so we will only
look at the cases where it seems likely that the solutions to the equations
will be easy to find.

Suppose we have three elliptic curves $E_1$, $E_2$, and $E_3$ over a field $k$
with each $E_i$ defined by an equation $y^2 = x(x^2 + A_i x + B_i)$
and such that the product $\Delta_1\Delta_2\Delta_3$ is a square in $k$.
For each $i$ let $Q_i$ be the point $(0,0)$ on $E_i$ and let $P_i$ be
some other $2$-torsion point on $E_i$, corresponding as in Section~\ref{genus3}
to a square root $d_i$ of $\Delta_i$.
We noted in the proof of Proposition~\ref{hyperelliptic3} that
the condition that the twisting factor be $0$ is equivalent to
the condition that 
\begin{equation}
\label{lambdaeq}
\lambda_3 = \lambda_1\lambda_2\pm\sqrt{(\lambda_1^2-1)(\lambda_2^2-1)},
\end{equation}
where $\lambda_i = A_i / d_i$. 
We can rewrite this equation in the equally useful form
\begin{equation}
\label{Aoverd}
\frac{A_3}{d_3} = \frac{A_1A_2 \pm4 \sqrt{B_1B_2}}{d_1d_2}.
\end{equation}
The fact that these equations hold precisely when $T=0$
will be the basis of all of our constructions of
genus-$3$ hyperelliptic curves with large torsion subgroups.

\subsection{Building hyperelliptic Jacobians with $E_3$ of type $(2,2)$}

Suppose $E_3 = F_{2,2}^t$, so that $A_3 = -t-1$ and $B_3 = t$. We have 
$\Delta_3 = (t-1)^2$, so let us choose $d_3$ to be $1-t$.  If $E_1$ and $E_2$
are any elliptic curves over $\Q$ such that $B_1 B_2$ and $\Delta_1\Delta_2$ are
both squares, say $\Delta_1\Delta_2 = r^2$ and $B_1B_2=s^2$,
then equation~(\ref{Aoverd})
becomes
$$\frac{t+1}{t-1} = \frac{A_1A_2 \pm 4 s}{r}.$$
This will have a rational nonzero solution for $t$ as long as the right-hand side is
neither $1$ nor $-1$.  Thus, we need only search for pairs $(E_1, E_2)$ such that
$B_1B_2$ and $\Delta_1\Delta_2$ are both squares.

Take $E_1 = F_{10}^2$ and $E_2 = F_6^u$.  Then $B_1 = -2^9$ and
$\Delta_1 = 11\cdot 3^5$, and up to squares in $\Q(u)$ we have 
\begin{align*}
         B_1 B_2 &=  2u \\
\Delta_1\Delta_2 &=  33(9u+1)/(u+1).
\end{align*}
Thus we would like to find rational solutions to the pair of equations
$u= 2v^2$ and $(9u+1)/(u+1) = 33w^2$.  Solving the second equation for
$u$ gives $u = (33w^2-1)/(9-33w^2)$, and inserting this in the
first equation and setting $x = 2(9-33w^2)v$ gives us
$$x^2 = 2(33w^2 - 1)(9-33w^2).$$
This curve of genus~$1$ has a rational point $(w,x) = (1/3,16/3)$, and a calculation
shows that it is birational with the elliptic curve
$$y^2 = x(x+66)(x-198).$$
This elliptic curve has rank $2$; its group of rational points is
generated by its $2$-torsion and the points $(-44,484)$ and $(-2, 160)$.
Suppose $C$ is the curve associated to one of these rational points
via  Proposition~\ref{hyperelliptic3} and our choice of the curves $E_i$.
We leave it to the reader to show that the image of the rational
torsion of $E_1\times E_2\times E_3$ in the Jacobian of $C$ 
is a group of the form $\Z/2\Z \times \Z/30\Z$.
Thus we have a family of hyperelliptic curves of genus~$3$,
parameterized by the points on a positive rank elliptic curve, whose
Jacobians contain a rational subgroup of the form $\Z/2\Z \times \Z/30\Z$.

Suppose we take $E_1$ to be equal to $E_2$.  Then $B_1B_2$ and $\Delta_1\Delta_2$
are automatically squares.  In particular, if we take $E_1 = F_{10}^t$,
we find a $1$-parameter family of hyperelliptic genus-$3$ curves
whose Jacobians contain $\Z/10\Z\times\Z/10\Z$.

Suppose we take $E_1 = F_{8,2}^u$ and $E_2 = F_{8,2}^v$.
Then $\Delta_1$ and $\Delta_2$ are each squares in $\Q(u,v)$.  To find
values of $u$ and $v$ that make the product $B_1B_2$ a square,
we must find rational solutions to the equation
$$(u^2 - 2u - 1)(u^2 + 2u - 1) w^2 = (v^2 - 2v - 1)(v^2 + 2v -1).$$
This equation defines an elliptic surface over $\Q(u)$ if we take
the zero section to be $(v,w)=(u,1)$.
Then the section $(v,w)=(-u,1)$ has infinite order.
Thus we find a family of hyperelliptic genus-$3$ curves,
parameterized by the points on a positive rank elliptic surface,  whose Jacobians contain 
a rational subgroup of the form $\Z/2\Z \times \Z/8\Z \times \Z/8\Z$.

The groups we can obtain from other choices of $E_1$ and $E_2$ are subgroups of 
the groups we build in the next few sections.

\subsection{Building hyperelliptic Jacobians with $E_3$ of type $(2,4)$}

Suppose now we take $E_3 = F_{2,4}^t$, so that $A_3 = 2t^2+2$,
$B_3 = (t-1)^2(t+1)^2$, and $\Delta_3 = 16t^2$.  If we take
$d_3=4t$, then equation~(\ref{Aoverd}) becomes
$$\frac{t^2+1}{2t} = \frac{A_1A_2 \pm4 \sqrt{B_1B_2}}{\sqrt{\Delta_1\Delta_2}}.$$
Solving this quadratic equation for $t$, we find
$$ t = \frac{ (A_1\pm 2\sqrt{B_1}) (A_2\pm 2\sqrt{B_2})} {\sqrt{\Delta_1\Delta_2}}.$$
Thus, if $B_1$, $B_2$, and $\Delta_1\Delta_2$ are all squares in $\Q$, we can
find a specialization of $F_{2,4}^t$ that will give us a twisting factor of $0$.
(Note that $A_i \not= \pm 2 \sqrt{B_i}$ since $E_i$ is nonsingular,
so the bad value $t=0$ in $F_{2,4}^t$ is automatically avoided.)

Suppose we take $E_1 = F_{8}^u$ and $E_2 = F_{12}^v$.  We check from
Table~\ref{universal2} that $B_1$ and $B_2$ are squares in $\Q(u,v)$,
so all we must do is find values of $u$ and $v$ such that 
$\Delta_1\Delta_2$ is a square.  Finding such $u$ and $v$ reduces to finding
rational points on the surface $S$ defined by
$$ (2u^2 - 1) w^2 = (6v^2-6v+1)(2v^2-2v+1).$$
Let $Y$ be the genus-0 curve $2t^2-1=41z^2$, and let $E$ be the genus-1 curve
	$$41y^2 = (6v^2-6v+1)(2v^2-2v+1).$$
We have a rational map $Y \times E \rightarrow S$ over $\Q$
mapping $(t,z),(v,y)$ to $(u,v,w)=(t,v,y/z)$.
Since $Y$ has the rational point $(t,z)=(9/11,1/11)$, $Y$ is isomorphic
to $\PP^1$ over $\Q$.
Since $E$ has the rational point $(v,y)=(5,11)$, it is an elliptic curve
over $\Q$, and in fact it is isomorphic to $y^2=x^3-41x^2+1681x$.
Moreover $E$ has positive rank, since $x=81/121$ gives a point of
infinite order.
Hence $Y \times E$ is a split elliptic surface over $\PP^1_\Q$ of positive rank,
and the points on this surface parameterize a family of hyperelliptic
curves of genus~$3$ whose Jacobians contain groups isomorphic
to $\Z/2\Z \times \Z/2\Z \times \Z/2\Z \times \Z/24\Z$.

Next we take $E_1=F_{2,8}^u$ and $E_2 = F_{2,8}^v$.
We check that $B_1$, $B_2$, $\Delta_1$, and $\Delta_2$ are all
squares in $\Q(u,v)$, so every pair of rational values of $u$ and $v$ will
give us a rational values of $t$.  This gives us a $2$-parameter
family of hyperelliptic curves of genus~$3$ whose Jacobians contain groups isomorphic
to $\Z/2\Z \times \Z/2\Z \times \Z/2\Z \times \Z/4\Z \times \Z/8\Z$.

Finally, if we take $E_1 = F_{12}^u$ and $E_2 = F_{12}^v$, we see that 
$B_1$ and $B_2$ are squares.  The condition that $\Delta_1\Delta_2$ be a square
leads us to find rational solutions to the equation
$$(6u^2-6u+1)(2u^2-2u+1) w^2 = (6v^2-6v+1)(2v^2-2v+1).$$
This equation defines an elliptic surface over $\Q(u)$ if we take the
zero section to be $(v,w)=(u,1)$, and then the section $(v,w)=(-u+1,1)$ has infinite order.
Thus we find a positive rank elliptic surface whose points parameterize a
family of hyperelliptic curves of genus~$3$ whose Jacobians contain groups
isomorphic to $\Z/2\Z \times \Z/2\Z \times \Z/6\Z \times \Z/12\Z$.

\subsection{Building hyperelliptic Jacobians with $E_3$ of type $(4,2)$}

Building hyperelliptic Jacobians by taking $E_3 = F_{4,2}^t$ is much more
difficult than doing so by taking $E_3 = F_{2,4}^t$, but it is possible.
Here we have $A_3 = -t^2 - 6t - 1$ and we may take $d_3 = (t-1)^2$, so
equation~(\ref{Aoverd}) becomes
\begin{equation}
\label{hyper42eq}
\frac{-t^2 - 6t - 1}{t^2 - 2t + 1} = \frac{A_1A_2 \pm4 \sqrt{B_1B_2}}{d_1d_2}.
\end{equation}
Suppose we have chosen elliptic curves $E_1$ and $E_2$, and have calculated
the right hand side of equation~(\ref{hyper42eq}) to get a number $r$.
Solving for $t$, we find that we must have $(r+1)t^2 + (2-6r)t + (r+1) = 0$,
and for $t$ to be a rational number the discriminant of this quadratic must be
a square, which reduces to the condition that $2(1-r)$ be a square.

If we set $E_1$ and $E_2$ equal to some of our universal elliptic curves,
the condition that $2(1-r)$ be a square turns into an absolute mess
that the reader should be thankful we do not go into here.
However, in the special case where $E_1 = E_2$ and we take the
plus sign in equation~(\ref{hyper42eq}), we have $r = 2\lambda_1^2 - 1$,
so $2(1-r)$ is a square precisely when $1-\lambda_1^2$ is a square.
A little algebra shows that this will be the case when
$-B_1\Delta_1$ is a square.

Suppose we take $E_1 = F_{8}^t$.  Finding $t$ such that $-B_1\Delta_1$
is a square is equivalent to solving the equation $w^2 = 1 - 2t^2$.
This equation defines a rational curve, so we obtain a $1$-parameter 
family of hyperelliptic curves, and we compute that the rational torsion
subgroup of the Jacobian of each of these curves contains a group isomorphic to
$\Z/4\Z \times \Z/4\Z \times \Z/8\Z$.

We can take $E_1$ to be any of the other universal curves, but for most of
the choices it is not possible to have $-B_1\Delta_1$ be a square, and
for the others the groups we get are subgroups of groups that we have already
obtained.

\subsection{Gaining $2$-power torsion}
\label{gaininghyper}

We noted in Section~\ref{gaining} that the rational torsion
subgroup of a quotient of an abelian variety is sometimes larger than 
the image of the rational torsion of the original variety.  As in the genus-$2$
case, we can use this fact to increase the size of the torsion subgroups we can make. 
In order to do this, it will be useful to have some of the ideas used in Section~\ref{gaining}
spelled out in more detail.

Let $k$ be a field with separable closure $K$ and
suppose $E$ is an elliptic curve over $k$ with $\#E[2](k) = 4$.
Let $y^2 = f(x) = (x-x_S)(x-x_T)(x-x_U)$ be a model for $E$, and let $S$, $T$, and $U$ be
the $2$-torsion points on $E$ with $x$-coordinates $x_S$,  $x_T$,  and $x_U$, respectively.

\begin{lemma}
\label{galoisaction}
Let notation be as above, and
suppose $W$ is an element of $E(K)$ such that $2W = S$.
Then $W$ can be defined over the field $\ell = k(\sqrt{x_S-x_T}, \sqrt{x_S-x_U})$.
Furthermore, the action of an element $\sigma$ of $\Gal(K/k)$ on $W$ can be determined
by its action on $\sqrt{x_S-x_T}$ and $\sqrt{x_S-x_U}$ as follows:

\begin{enumerate}
\item[(a)] If $\sigma$ fixes neither $\sqrt{x_S-x_T}$ nor $\sqrt{x_S-x_U}$, then $W^\sigma-W = S$.
\item[(b)] If $\sigma$ fixes $\sqrt{x_S-x_T}$ but not $\sqrt{x_S-x_U}$, then $W^\sigma-W = T$.
\item[(c)] If $\sigma$ fixes $\sqrt{x_S-x_U}$ but not $\sqrt{x_S-x_T}$, then $W^\sigma-W = U$.
\end{enumerate}
\end{lemma}

\begin{proof}
Under our assumptions, the $k$-algebra $L$ of Section~\ref{gaining}
is isomorphic to $k\times k\times k$.
The lemma then follows from the fact that the isomorphism 
$$H^1(G_k,E[2]) \isom \ker\left(L^*/L^{*2} 
   \stackrel{\Norm}\longrightarrow k^*/k^{*2} \right)$$ from \cite{schaefer}
sends the image of $S$ in $H^1(G_k,E[2])$ to the class of the element
$$\big( (x_S-x_T)(x_S-x_U), (x_S-x_T), (x_S-x_U) \big).$$
\end{proof}

Using Lemma~\ref{galoisaction} we can write down families of elliptic curves
with specific rational torsion subgroups and extra $4$-torsion over an 
extension that is at most quadratic.  For instance, suppose we take $E$ to be
the universal curve $F_{2,4}^t$, with $x_S = -(1+t)^2$ and $x_T = 0$ and $x_U=-(1-t)^2$.
We see that $x_S-x_T = -(1+t)^2$ and $x_S-x_U = -4t$.  If we take $t=-s^2$ then
$x_S-x_U$ will be a square and $x_S-x_T$ will be $-1$, up to squares, so for
this choice of $t$ there will be a point $W$ of $E$ defined over $\ell=k(\sqrt{-1})$
such that $2W=S$ and such that $W^\sigma - W = U$ for every non-trivial $k$-automorphism
$\sigma$ of $\ell$.  With a little calculation we can find the coordinates for $W$.
Putting this together with the information we have about $F_{2,4}^t$ from
Tables~\ref{universal2}, \ref{maximal}, and~\ref{explicit2torsion2} gives us
the information summarized in Table~\ref{moretorsion24}.  We will refer to
the universal curve we have thus constructed as $F_{2,4a}^s$.

\begin{table}
\begin{center}
\begin{tabular}{|c|l|}								\hline
Coefficients of model			& $A = 2  (s^4 + 1)$ 				\buff\\	
$y^2 = x(x^2+Ax+B)$    			& $B = (s^2 + 1)^2  (s + 1)^2  (s - 1)^2$	\buff\\	
of universal curve    			& $\Delta = 16s^4$ 				\buff\\ 
\hline\hline
$x$-coordinates	of			& $x_S = -(s - 1)^2 (s + 1)^2$ 			\buff\\	
$2$-torsion points			& $x_T = 0$ 					\buff\\	
 $S$, $T$, $U$	 			& $x_U = -(s^2 + 1)^2$ 				\buff\\
\hline\hline
$x$- and $y$-coordinates of a		& $x_V = - (s^2 + 1) (s + 1) (s - 1)$ 		\buff\\
$4$-torsion point $V$ with $2V=T$	& $y_V =  2 (s^2 + 1) (s + 1) (s - 1)$ 		\buff\\ 
\hline\hline
$x$- and $y$-coordinates of a		& $x_W = - (s - 1)  (s + 1)  (s - i)^2$ 	\buff\\
$4$-torsion point $W$ with $2W = S$	& $y_W = -2  s  (s - 1)  (s + 1)  (s - i)^2$ 	\buff\\ 
\hline
\end{tabular}
\end{center}
\vspace{1ex}
\caption{Data for the universal elliptic curve $F_{2,4a}^s$.
Here $i$ denotes a square root of $-1$. If $\sigma$ is a non-trivial
element of $\Gal(k(i)/k)$, then $W^\sigma-W = U$.}
\label{moretorsion24}
\end{table}

We can use the curve $F_{2,4a}^s$ to build Jacobians.  Suppose we take
$E_1 = E_2$, $E_3 = F_{2,4a}^s$, $d_1 = d_2$, and $d_3 = 4s^2$.
Then equation~(\ref{lambdaeq}) gives us two possibilities:
either $\lambda_3 = 1$,
or $\lambda_3 = 2\lambda_1^2 - 1$.
The former is impossible, because in that case $x(x^2+A_3x+B_3)$
would have a double root at~0.
The latter becomes
$$\frac{s^4+1}{2s^2} = 2\lambda_1^2 - 1,$$
which can be solved to obtain
$s = \pm \lambda_1 \pm\sqrt{\lambda_1^2 - 1}$, or
$$s = \frac{\pm A_1 \pm 2\sqrt{B_1}}{\sqrt{\Delta_1}}.$$
So suppose we take $E_1=E_2 = F_{2,8}^u$ over $k = \Q(u)$.
For this curve the numbers $B_1$ and $\Delta_1$
are both squares in $k$, so we can set $s$ to be an element of $k$
that makes the twisting factor equal to zero.
Thus we find a hyperelliptic curve $C$ over $k=\Q(u)$ whose Jacobian 
is $(2,2,2)$-isogenous to $E_1\times E_1\times E_3$.
Given the choice of points $Q_i$ implicit
in the above expressions,  it is easy to calculate that the image of the known
$k$-rational torsion of $E_1\times E_1\times E_3$ in the Jacobian $J$ of $C$
is a group of the form $\Z/2\Z \times \Z/2\Z \times \Z/2\Z \times \Z/4\Z \times \Z/8\Z$.
Now we will show that in fact $J(k)$ contains a torsion group larger than this.

Let $P_1$ be the point on $E_1$ with $x = -16u^4$ and let $R_1$ be
the point with $x = - (u-1)^4 (u+1)^4$  (see Table~\ref{explicit2torsion2}).
If we apply Lemma~\ref{galoisaction} to the curve $E_1$, with 
$S=P_1$ and $T=Q_1$ and $U=R_1$, we find that there is a point $W_1\in E_1(K)$
defined over an (at worst) biquadratic extension $\ell$ of $k$ with $2W_1 = P_1$.
Note that $x_S - x_T = -16u^4$ differs from $-1$ by a square, so $k(i)\subseteq\ell$
and the action of a $\sigma\in\Gal(\ell/k)$ on $\sqrt{x_{P_1}-x_{Q_1}}$ is
the same as its action on $i$.

Let $P_3$ be the point on $E_3$ with $x = -(s^2+1)^2$ (see Table~\ref{moretorsion24})
and let $W_3$ be the point labeled $W$ in Table~\ref{moretorsion24}.
Let us consider how an element $\sigma$ of $\Gal(\ell/k)$ acts upon
the element $(W_1, W_1, W_3)$ of $(E_1\times E_1\times E_3)(\ell)$.
If $\sigma$ is not the identity and yet fixes $i$, then
$(W_1, W_1, W_3)^\sigma - (W_1,W_1,W_3) = (Q_1, Q_1, 0)$.
If $\sigma$ does not fix $i$, then $W_3^\sigma - W_3 = P_3$ and we see
that $(W_1, W_1, W_3)^\sigma - (W_1,W_1,W_3)$ is either $(P_1,P_1,P_3)$ or
$(R_1,R_1,P_3)$.  Thus we see that $(W_1, W_1, W_3)^\sigma - (W_1,W_1,W_3)$
is an element of the kernel $G$ of $E_1\times E_1\times E_3 \ra J$ for
every $\sigma$, so the image $Z$ of $(W_1, W_1, W_3)$ in $J$ is defined over $k$.
Since $2(W_1,W_1,W_3) = (P_1,P_1,P_3)$ is in $G$, the point $Z$ is a $2$-torsion point.
Finally, we note that $(W_1, W_1, W_3)$ does not differ by an element of $G$ from
any of the known rational torsion points of $E_1\times E_1\times E_3$, so $Z$
is not in the image of the known rational torsion of  $E_1\times E_1\times E_3$.
Thus, the $k$-rational torsion of the Jacobian of $C$ contains a group
isomorphic to $\Z/2\Z \times \Z/2\Z \times \Z/2\Z \times \Z/2\Z \times \Z/4\Z \times \Z/8\Z$.
Since $k$ is a rational function field over $\Q$, we get a $1$-parameter family
of hyperelliptic curves of genus~$3$ having this group in their
rational torsion subgroups.

\subsection{Building Jacobians of plane quartics --- introduction}

Now we turn our attention to the task of building plane quartics whose
Jacobians have large rational torsion subgroups.  
Most of the families we will construct will be produced by
fixing the torsion structure of the curve $E_3$ and analyzing
the twisting factor as a function of the coefficients of $E_1$ and $E_2$,
but one interesting class of examples will arise by setting $E_1=E_2=E_3$.
Instead of trying to get the twisting factor to be zero, as we did in
the last few sections, we will try to get the twisting factor to be a square,
so that (by the final statement of Proposition~\ref{quartic3}) 
the product of the $E_i$ will be isogenous over $\Q$ to the
Jacobian of a curve.
It will be convenient to use the second expression for the twisting factor,
namely
\begin{equation}
\label{lambdaeq2}
T = d_1 d_2 d_3(\lambda_1^2 + \lambda_2^2 + \lambda_3^2 - 2\lambda_1\lambda_2\lambda_3 -1),
\end{equation}
where $\lambda_i = A_i /d_i$.

As in the hyperelliptic case, we will not examine all possible combinations of torsion 
structures on the curves $E_1$, $E_2$, and $E_3$ here because of the
complexity of the equations that arise.

\subsection{Building Jacobians of plane quartics with $E_3$ of type $(4,2)$}
\label{building42}

In this section we will take $E_3 = F_{4,2}^t$, so that $A_3 = -t^2 - 6t - 1$
and $\Delta_3 = (t-1)^4$.  We will take $d_3 = (t-1)^2$.

Suppose we take $E_1 = E_2$ and $d_1=d_2$.  Then $\lambda_1 = \lambda_2$, and
the twisting factor is 
\begin{align*}
T &= d_1^2 d_3 (2\lambda_1^2 + \lambda_3^2 - 2\lambda_1^2\lambda_3 - 1)\\
  &= \Delta_1 d_3 (\lambda_3 - 1)(\lambda_3 - 2\lambda_1^2 + 1) \\
  &= \Delta_1 (\lambda_3 - 1)(A_3 - (2\lambda_1^2 - 1)d_3) \\
  &= (\lambda_3 - 1) (\Delta_1 A_3 - (2A_1^2 - \Delta_1) d_3).
\end{align*}
A quick calculation shows that $\lambda_3 - 1 = -2 (t+1)^2 / (t-1)^2$,
so up to squares in $\Q(t)$ the twisting factor is
$$4A_1^2 d_3  -2\Delta_1(A_3+d_3) = 4(t-1)^2 A_1^2 + 16 t \Delta_1.$$
To get the twisting factor to be a square, we need to find rational solutions
to the equation $$w^2 = 4(t-1)^2 A_1^2 + 16 t \Delta_1.$$
For fixed $A_1$ and $\Delta_1$, this last equation defines a curve of genus~$0$
in the $(t,w)$-plane, and since it has a rational point
(namely $(t,w)=(0,2A_1)$) it is isomorphic to $\PP^1$.  We can parameterize
the curve by setting
\begin{align*}
	t &= (z+4B_1)(z-\Delta_1)/(A_1^2 z), \\
	w &= 2(z^2 + 4B_1\Delta_1)/(A_1z).
\end{align*}

Suppose in particular we take $E_1 = F_{10}^u$.  Then we get a $2$-parameter family
(the parameters being $u$ and $z$) of plane quartics whose Jacobians are
isogenous to $E_1\times E_1\times E_3$, and a simple computation with abelian
groups shows that these Jacobians have a rational subgroup isomorphic to 
$\Z/10\Z\times\Z/20\Z$.

If we take $E_1 = F_{2,8}^u$ we get a $2$-parameter family of plane quartics whose
Jacobians contain a rational subgroup isomorphic to 
$\Z/2\Z\times\Z/4\Z\times\Z/4\Z\times\Z/8\Z$.
If we take $E_1 = F_{8,2}^u$ we get a $2$-parameter family of plane quartics whose
Jacobians contain a rational subgroup isomorphic to 
$\Z/2\Z \times\Z/2\Z \times \Z/8\Z \times \Z/8\Z$.
And if we take $E_1 = F_{12}^u$ we  get a $2$-parameter family of plane quartics whose
Jacobians contain a rational subgroup isomorphic to 
$\Z/2\Z \times\Z/12\Z \times \Z/12\Z$.

Without the assumption that $E_1 = E_2$ it is not as easy to make the twisting
factor a square.
But suppose we take $E_1 = F_{10}^{-1/2}$ and $E_2 = F_8^{1/2}$.
Then $\lambda_1^2 = -625/2048$ and
$\lambda_2^2 = -49/32$, and we can choose $d_1$ and $d_2$ so that
$d_1d_2 = 4$ and $\lambda_1\lambda_2 = 175/256$.
Then the twisting factor is
\begin{align*}
T & = 4 d_3 \left( -\frac{625}{2048} - \frac{49}{32} + \frac{A_3^2}{\Delta_3}
       - \frac{175}{128} \frac{A_3}{d_3} - 1\right) \\
  & = \frac{-1922t^4 + 118024t^3 + 29940t^2 + 118024t - 1922}{2^{10}(t-1)^2},
\end{align*}
so in order to make the twisting factor a square we must find rational solutions
to the equation
$$w^2 = -1922t^4 + 118024t^3 + 29940t^2 + 118024t - 1922.$$
This last equation defines a curve of genus~$1$, and it has a rational
point, namely $(t, w) = (1, 512)$.  A calculation then shows that the curve
is birational with the elliptic curve defined by 
$y^2 = x^3 + 2565x -15606$.  This happens to be the curve 528A2 in \cite{cremona},
which has rank $1$.  (The point $(33,324)$ is of infinite order.)
We see that there is a family of plane quartics, parameterized by
the points on a positive rank elliptic curve, 
whose Jacobians contain a rational subgroup isomorphic to
$\Z/4\Z\times\Z/40\Z$.

Now let us try taking $E_1 = F_{8,2}^2$ and $E_2 = F_{6,2}^{1/4}$.
If we take $d_1 = 256$ and $d_2 = 4$ then we have 
$\lambda_1 = 47/128$ and $\lambda_2 = 863/512$.
The twisting factor is
\begin{align*}
T & = 1024 d_3 \left( \frac{ 2209}{16384} + \frac{744769}{262144} + \frac{A_3^2}{\Delta_3}
       - \frac{40561}{32768} \frac{A_3}{d_3} - 1\right) \\
  & = \frac{1104601t^4 + 2371804t^3 + 9824406t^2 + 2371804t + 1104601}{2^8(t-1)^2},
\end{align*}
so in order to make the twisting factor a square
we must find rational solutions to the equation
$$w^2 = 1104601t^4 + 2371804t^3 + 9824406t^2 + 2371804t + 1104601.$$
The genus-$1$ curve defined by this equation has a rational point ---
namely, $(t,w)=(1,4096)$ --- so it is an elliptic curve.
A calculation shows that it is birational with the elliptic curve
$y^2 = x^3 -151563x +  10810438$.
Cremona's rank-finding program calculates that
this elliptic curve has rank $2$, and provides the two independent rational
points $(59,1440)$ and $(-157,5544)$.
Thus we have a positive rank elliptic curve whose points parameterize a
family of plane quartics whose Jacobians contain a rational subgroup isomorphic to
$\Z/2\Z\times\Z/4\Z\times\Z/24\Z$.

Finally, let us take $E_1 = F_{10}^{-1/3}$ and $E_2 = F_{12}^{1/3}$.
Then $\lambda_1^2 = -485809/759375$ and
$\lambda_2^2 = -3721/375$, and we can choose $d_1$ and $d_2$ so that
$d_1d_2 =  625/59049$ and $\lambda_1\lambda_2 = -42517/16875$.
The twisting factor is
\begin{align*}
T & =  \frac{625}{59049}d_3 \left( -\frac{485809}{759375} - \frac{3721}{375} 
       + \frac{A_3^2}{\Delta_3}
       + \frac{85034}{16875} \frac{A_3}{d_3} - 1\right) \\
  & = \frac{-177710460t^4 + 433908240t^3 + 216604440t^2 + 433908240t - 177710460}
            {3^{16} 5^2 (t-1)^2},
\end{align*}
so in order to make the twisting factor a square we must find rational solutions
to the equation
$$w^2 = -177710460t^4 + 433908240t^3 + 216604440t^2 + 433908240t - 177710460.$$
This genus-$1$ curve has
rational points --- for instance, $(t,w)=(1,27000)$ ---
so it can be made into an elliptic curve.  A calculation shows that the
elliptic curve we get is isomorphic to the curve $y^2 = x^3 + 213x -30566$.
Cremona's rank-finding program says that this curve has rank $1$, and
gives the non-torsion point $(53,360)$.
Thus we have found a positive rank elliptic curve whose points
parameterize a family of plane
quartics whose Jacobians contain a rational subgroup isomorphic to
$\Z/4\Z\times\Z/60\Z$.

\subsection{Building Jacobians of plane quartics with $E_3$ of type $(2,4)$}

In this section we will take $E_3 = F_{2,4}^t$, so that $A_3 = 2t^2+2$
and $\Delta_3 = 16t^2$.  We will take $d_3 = 4t$.

As we saw in the preceding section, if we take $E_1 = E_2$ and $d_1 = d_2$ then
the twisting factor is
$$T = (\lambda_3 - 1) (\Delta_1 A_3 - (2A_1^2 - \Delta_1) d_3).$$
For our choice of $E_3$ we have $\lambda_3 - 1 = (t-1)^2 / (2t)$,
so up to squares in $\Q(t)$ the twisting factor is
$$(2t)(\Delta_1 (2t^2+2) - (2A_1^2 -\Delta_1)(4t)) =
  4t (\Delta_1 (t+1)^2 - 4 A_1^2 t).$$
Thus, we would like to find rational solutions to the 
equation
\begin{equation}
\label{24eq}
w^2 = 4t (\Delta_1 (t+1)^2 - 4 A_1^2 t).
\end{equation}

Suppose we take $E_1 = F_{8,2}^2$, so that $A_1 = 2\cdot 47$ and $\Delta_1 = 2^{16}$.
Equation~(\ref{24eq}) becomes
$w^2 = 4t (2^{16}(t+1)^2 - 2^4 47^2 t),$ and by setting $s = 4096t$ and $z=512w$ we get
$$z^2 = s^3+ 5983 s^2 + 16777216s.$$
A search for points on this curve using Cremona's programs comes up with
the non-torsion point $(s,z)=(5929/64,20520885/512)$.
Thus we find a family of plane quartics, parameterized by the
points on a positive rank elliptic curve,  whose Jacobians contain
a rational subgroup isomorphic to
$\Z/4\Z\times \Z/8\Z\times\Z/8\Z$.

Other choices for $E_1$ lead to groups we have already constructed,
most of them in the sections on hyperelliptic curves.

\subsection{Building Jacobians of plane quartics with $E_1=E_2=E_3$}

Suppose we take $E_1=E_2=E_3$ and $d_1 = d_2 = d_3$.
Since $\Delta_1\Delta_2\Delta_3$ is supposed
to be a square, we see that $\Delta_1$ must be a square and $d_1$ must an 
element of the base field.  If we write $\lambda$ and  $d$ for
$\lambda_1$ and $d_1$, we find that the twisting factor is
$$T = d^3 (3\lambda^2 -2\lambda^3 -1) =- d^3(\lambda-1)^2(2\lambda+1).$$
Up to squares, then, the twisting factor is
$$-d(2\lambda+1) = -(2A+d),$$
where we write $A$ for $A_1$.

Suppose we take $E_1 = F_{2,6}^t$.  Then $A = -2t^4+12t^2+6$ and we can take
$d = 16t$.  The twisting factor, up to squares, is
$t^4 - 6t^2 - 4t -3$, so we would like to find rational solutions to
$$w^2 = t^4 - 6t^2 - 4t -3.$$
The desingularization of this genus-$1$ curve has rational points at infinity,
so it is an elliptic curve.
A calculation shows that it is isomorphic as an elliptic curve to
$y^2 = x^3 - 48$. This is the curve 243A1 in \cite{cremona};
it has rank~$1$, and its group of rational points is generated by the
point~$(4,4)$.  
Thus we find a positive rank elliptic curve whose points
parameterize a family of plane quartics whose
Jacobians contain a rational subgroup isomorphic to
$\Z/6\Z \times \Z/6\Z \times \Z/6\Z$.

Taking $E_1$ to be $F_{2,4}^t$ or $F_{4,2}^t$ gives us rational subgroups
that we can get in other ways, and taking $E_1$ to be $E_{2,8}^t$ or
$E_{8,2}^t$ leads to equations with no rational solutions.

\subsection{Gaining $2$-power torsion}

In this section we will show how to further specialize two of the families we
wrote down in the last few sections to obtain larger rational
torsion subgroups.  For the first example, we will need to write down the
universal curve $F_{2,4a}^s$ in a slightly different form.  All we would
like to do is make a change of variables in the equation for $F_{2,4a}^s$
so that the point labeled $S$ in Table~\ref{moretorsion24} will
have $x$-coordinate $0$.  After translating $x$ by the proper amount to
do this, we obtain a curve that we will call $F_{4,2a}^s$.  All the information
we will need about this curve is listed in Table~\ref{moretorsion42}. 
Note that $F_{4,2a}^s$ is a specialization of $F_{4,2}^t$, as the notation suggests.
In fact, $F_{4,2a}^s  = F_{4,2}^{-s^2}$.

\begin{table}
\begin{center}
\begin{tabular}{|c|l|}								\hline
Coefficients of model			& $A = -(s^2 - 2s - 1)(s^2 + 2s - 1)$		\buff\\	
$y^2 = x(x^2+Ax+B)$    			& $B = -4 s^2 (s - 1)^2 (s + 1)^2$		\buff\\	
of universal curve    			& $\Delta = (s^2 + 1)^4$ 			\buff\\ 
\hline\hline
$x$-coordinates	of			& $x_S = 0$		 			\buff\\	
$2$-torsion points			& $x_T = (s + 1)^2 (s - 1)^2$ 			\buff\\	
 $S$, $T$, $U$	 			& $x_U = -4s^2$ 				\buff\\
\hline\hline
$x$- and $y$-coordinates of a		& $x_V = - 2 (s + 1) (s - 1)$ 			\buff\\
$4$-torsion point $V$ with $2V=T$	& $y_V =  2 (s^2 + 1) (s + 1) (s - 1)$ 		\buff\\ 
\hline\hline
$x$- and $y$-coordinates of a		& $x_W = 2 i s  (s + 1)  (s - 1)$	 	\buff\\
$4$-torsion point $W$ with $2W = S$	& $y_W = -2  s  (s + 1)  (s - 1)  (s - i)^2$ 	\buff\\ 
\hline
\end{tabular}
\end{center}
\vspace{1ex}
\caption{Data for the universal elliptic curve $F_{4,2a}^s$.
Here $i$ denotes a square root of $-1$. If $\sigma$ is a non-trivial
element of $\Gal(k(i)/k)$, then $W^\sigma-W = U$.}
\label{moretorsion42}
\end{table}

Now suppose we try to build a plane quartic by taking $E_1 = E_2 = F_{2,8}^u$
and $E_3 = F_{4,2a}^s$.  For notational convenience, we denote by $t$ the number $-s^2$,
so that $E_3 = F_{4,2}^t$ as well.  
Let $P_1$ be the $2$-torsion point on $E_1$ with $x=-16u^4$
and let $P_3$ be the point labeled $U$ in Table~\ref{moretorsion42},
so that $d_3 = (t-1)^2 = (s^2+1)^2$.

As we noted in Section~\ref{building42}, the twisting factor in this situation is
equal to
$$4(t-1)^2 A_1^2 + 16 t \Delta_1 = 4(s^2+1)^2 A_1^2 - 16s^2 \Delta_1,$$
up to squares, so we would like to find solutions to 
$$w^2 = (s^2+1)^2 A_1^2 - 4s^2 \Delta_1.$$
For fixed $A_1$ and $\Delta_1$, this is a genus-$1$ curve in the $(s,w)$-plane
whose desingularization has rational points at infinity, so it is isomorphic
to its Jacobian, which (according to the formulas in
Example 3.7 (pp.~293--294) of \cite{silverman}) is given by
$$y^2 = x (x+\Delta_1) (x - 4B_1).$$
If we take $u=2$ then this curve is given by
$y^2 = x (x+30625) (x-82944)$
and has a non-torsion point $(-21600,4514400)$.  
Thus, for this choice of $u$ there are infinitely many values of $s\in \Q$
that make the twisting factor a square.

Let $C$ be the plane quartic associated to one of these choices.
Let $A$ be the abelian variety $E_1\times E_1\times E_3$ and let $G$ be the
kernel of the homomorphism $\psi$ from $A$ to the Jacobian $J$ of $C$.
The image under $\psi$ of the known rational torsion of $A$ is a rational
torsion subgroup of $J$ isomorphic to 
$\Z/2\Z \times  \Z/4\Z\times \Z/4\Z \times \Z/8\Z$.  
Note that the number of independent $(\Z/4\Z)$-factors contained in this group
is as large as is allowed by the 
restrictions imposed by the Galois-equivariance of the Weil pairing.
However, there is more rational torsion on $J$ than just this.

Let $R_1$ be the $2$-torsion point on $E_1$ with $x = -(u-1)^4 (u+1)^4$
(see Table~\ref{explicit2torsion2}).
Applying Lemma~\ref{galoisaction} to the curve $E_1$, taking 
$S=P_1$ and $T=Q_1$ and $U=R_1$, we find that there is a point $W_1\in E_1(K)$
defined over a Galois extension $\ell$ of $\Q$ with $2W_1 = P_1$.
Since $x_S-x_T = -16u^4$ we see that $\ell$ contains $\Q(i)$;
furthermore, if $\sigma$ is any automorphism of $\ell$ that fixes $i$ then
$W_1^\sigma - W_1$ is either $0$ or $Q_1$, while if $\sigma$ does not fix $i$
then $W_1^\sigma - W_1$ is either $P_1$ or $R_1$.

Let $W_3$ be the point on $E_3$ labeled $W$ in Table~\ref{moretorsion42}.
The point $W_3$ is defined over $\Q(i)$, and if $\sigma$ is the non-trivial
automorphism of $\Q(i)$ then $W_3^\sigma - W_3=P_3$.
Knowing the Galois action on $W_1$ and $W_3$, one can check quite simply
that $(W_1, W_1, W_3)^\sigma - (W_1, W_1, W_3)$ is an element of $G$
for every $\sigma \in\Gal(K/\Q)$.  Thus the image $Z$ of $(W_1, W_1, W_3)$ in
$J$ is a rational point.  It is certainly not one of the points we
already calculated because $(W_1, W_1, W_3)$ does not differ from one of 
the known rational points of $A$ by an element of $G$.
However, $2Z$ is one of the points we previously calculated,
and indeed $2Z$ is a $2$-torsion point.
As remarked before, the torsion subgroup cannot contain $(\Z/4\Z)^n$
for $n>3$, so the only possibility is that $J$ contains a rational
torsion subgroup isomorphic to
$\Z/2\Z \times  \Z/2\Z\times \Z/4\Z \times \Z/4\Z \times \Z/8\Z$.  

We end by trying to add $2$-power torsion to the family of curves we obtained
in the preceding section by taking $E_1 = E_2 = E_3 = F_{2,6}^t$.
Since the three elliptic curves are isomorphic to one another,
we drop the subscripts.
Recall that we chose $d = 16t$.  This corresponds to choosing $P$
to be the $2$-torsion point on $E$ with $x$-coordinate equal to
$(t-3)(t+1)^3$.  Let $R$ be the $2$-torsion point on $E$ with
$x$-coordinate $(t+3)(t-1)^3$.  Of course, $Q$ is the point $(0,0)$ on $E$.

Suppose we can find a value of $t\in\Q$ that makes the twisting factor a square
and that equals $(3+s^2)/(1-s^2)$ for some $s\in\Q$.
Then the $x$-coordinate of $P$ will be a square, 
and if we apply Lemma~\ref{galoisaction} to $E$, taking
$S = P$ and $T=Q$ and $U = R$, we find that there will be point $W$ on $E$
defined over a quadratic extension $\ell$ of $\Q$ such that $2W = P$ and
$W^\sigma - W = Q$ for the non-trivial element of $\Gal(\ell/\Q)$.
Let $C$ be the plane quartic corresponding to this hypothetical $t$,
and let $\psi\colon E\times E\times E \ra J$ be the map to the
Jacobian of $C$ with kernel $G$ generated by $(0,Q,Q)$, $(Q,0,Q)$, and $(P,P,P)$.
The image under $\psi$ of the rational torsion of $E\times E\times E$ is
a group isomorphic to $\Z/6\Z \times \Z/6\Z \times \Z/6\Z$, but there
is still more rational torsion on $J$.  For consider the element
$Y_1 = (W,W,0)$ of $E\times E\times E$:  since $Y_1^\sigma - Y_1 = (Q,Q,0)$ is
in $G$, we see the image of $Y_1$ in $J$ is a rational point $Z_1$.
Similarly, the image $Z_2$ of $Y_2 = (W,0,W)$ is a rational point.  Since neither
$2Y_1$ nor $2Y_2$ is in $G$, while both $4Y_1$ and $4Y_2$ are zero,
we see that $Z_1$ and $Z_2$ are $4$-torsion points.  Moreover, they
are independent $4$-torsion points, because neither $2(Y_1+Y_2)$ nor $2(Y_1+3Y_2)$
is in $G$.
Thus $J$ contains a rational torsion subgroup isomorphic to 
$\Z/6\Z \times \Z/12\Z \times \Z/12\Z$.

To find such a curve $C$, we must find a value of $s$ such that $t = (3+s^2)/(1-s^2)$
will make the twisting factor a square.  From the preceding section, we
know that the twisting factor will be a square if there is a $w$ such that
$$w^2 = t^4 - 6t^2 - 4t -3.$$
Inserting our formula for $t$ into this equation and clearing denominators shows that we
want to find solutions to 
$$y^2 = -s^8 + 6s^4 + 56s^2 + 3,$$
with $s\neq \pm1$.
Amazingly enough, there are such solutions to this equation:
we can take $(s,y)=(\pm 1/5,1432/625)$.  These solutions give $t=19/6$.
Inserting the corresponding values for $A$, $B$, and $\Delta$ into the
formulas of Proposition~\ref{quartic3} leads us to the plane quartic $C$
defined by
$$15625(X^4 + Y^4 + Z^4) - 96914(X^2 Y^2 + X^2 Z^2 + Y^2 Z^2) = 0.$$
The Jacobian of this curve contains a rational torsion subgroup of order
$6\cdot12\cdot12 = 864$.  Both $E$ and $C$ have good reduction modulo $7$;
since the reduction $E'$ of $E$ modulo $7$ has $12$ points ($12$ being the 
only multiple of $12$ lying within the Weil bounds), and
since the reduction $J'$ of the Jacobian of $C$ is isogenous to $E'\times E'\times E'$,
we see that $J'$ has $12^3 = 1728$ points.  Therefore the rational torsion
subgroup of $J$ is isomorphic either to $\Z/6\Z\times \Z/12\Z\times \Z/12\Z$, 
to $\Z/12\Z\times \Z/12\Z\times \Z/12\Z$, to 
$\Z/2\Z\times \Z/6\Z\times \Z/12\Z \times\Z/12\Z$,
or to $\Z/6\Z\times \Z/12\Z \times\Z/24\Z$.
The second and third possibilities can be ruled out by looking at the 
action of Galois on the $4$-torsion of $E\times E\times E$.
The fourth can also be ruled out: the dual isogeny
$\widehat{\Psi^\ast}: J \rightarrow E \times E \times E$ would take
a $24$-torsion point on $J$ to a point of order at least~$12$ on
$E \times E \times E$,
since $\Psi^\ast \widehat{\Psi^\ast}$ is multiplication-by-2 on $J$,
but the torsion subgroup of $E \times E \times E$ is of exponent only~$6$.
Thus, we have found a single plane quartic $C$ such that the
rational torsion subgroup of $\Jac C$ is isomorphic to
$\Z/6\Z\times \Z/12\Z\times \Z/12\Z$.


\begin{thebibliography}{99}

\bibitem{baker}
{\sc H. F. Baker},
{\em Abelian functions: Abel's theorem and the allied theory of
theta functions},
Cambridge Univ. Press, 1897 (reprinted 1995).

\bibitem{brumerkramer}
{\sc A.\ Brumer and K.\ Kramer},
The rank of elliptic curves,
{\em Duke Math. J.} {\bf 44} (1977), no.~4, 715--743.

\bibitem{casselsgenus2} {\sc J.\ W.\ S.\ Cassels},
The Mordell-Weil group of curves of genus~$2$, in: M. Artin, J. Tate (eds.),
{\em Arithmetic and Geometry I}, Birkh\"{a}user, Boston, (1983), 27--60.

\bibitem{casselsquartic} {\sc J.\ W.\ S.\ Cassels},
The arithmetic of certain quartic curves,
{\em Proc. Roy. Soc. Edinburgh Sect. A} {\bf 100} (1985), nos.~3--4, 201--218.

\bibitem{cremona} {\sc J.\ Cremona},
{\em Algorithms for modular elliptic curves}, Cambridge Univ.\ Press, 1992.

\bibitem{delignerapoport}
{\sc P. Deligne and M. Rapoport},
Les sch\'emas de modules de courbes elliptiques,
in:
A. Dold and B. Eckmann (eds.),
{\em Modular Functions of One Variable II},
Lecture Notes in Math. {\bf 349},
Springer-Verlag, New York, (1973), 143--316.

\bibitem{fermigier} 
{\sc S. Fermigier}, 
Exemples de courbes elliptiques de grand rang sur $\Q(t)$ et sur $\Q$
poss\'edant des points d'ordre $2$,
{\em C. R. Acad. Sci. Paris S\'er. I Math.} {\bf 322} (1996), no.~10, 949--967.

\bibitem{flynnlarge}
{\sc E. V. Flynn},
Large rational torsion on abelian varieties, 
{\em J. Number Theory} {\bf 36} (1990), no.~3, 257--265.

\bibitem{flynnsequences}
{\sc E. V. Flynn},
Sequences of rational torsions on abelian varieties,
{\em Invent. Math.} {\bf 106} (1991), no.~2, 433--442.

\bibitem{frey}
{\sc G. Frey},
On elliptic curves with isomorphic torsion structures and
corresponding curves of genus~$2$,
in:
J. Coates and S. T. Yau (eds.),
{\em Elliptic Curves, Modular Forms, and Fermat's Last Theorem},
International Press, Boston, (1995), 79--98.

\bibitem{freykani}
{\sc G. Frey and E. Kani},
Curves of genus 2 covering elliptic curves and an arithmetical application,
in:
G. van der Geer, F. Oort, and J. Steenbrink (eds.),
{\em Arithmetic Algebraic Geometry},
Progress in Math. {\bf 89},
Birkh\"auser, Boston, (1991), 153--176.

\bibitem{hasegawa}
{\sc Y.\ Hasegawa},
Table of quotient curves of modular curves $X_0(N)$ with genus~$2$,
{\em Proc. Japan Acad. Ser. A Math. Sci.} {\bf 71} (1995), no.~10, 235--239.

\bibitem{howe}
{\sc E. W. Howe},
Constructing distinct curves with isomorphic Jacobians in characteristic zero,
{\em Internat. Math. Res. Notes} {\bf 1995}, no.~4, 173--180.

\bibitem{howeleprevostpoonen}
{\sc E. W. Howe, F. Lepr\'evost, and B. Poonen},
Sous-groupes de torsion d'ordres \'elev\'es de Jacobiennes d\'ecomposables de 
courbes de genre $2$,
{\em C. R. Acad. Sci. Paris S\'er. I Math.} {\bf 333} (1996), no.~9, 1031--1034.

\bibitem{kani-number}
{\sc E. Kani},
The number of curves of genus two with elliptic differentials,
{\em J. Reine Angew. Math.} {\bf 485} (1997), 93--121.

\bibitem{kani-existence}
{\sc E. Kani},
The existence of curves of genus two with elliptic differentials,
{\em J. Number Theory} {\bf 64} (1997), no.~1, 130--161.

\bibitem{katzmazur}
{\sc N. M. Katz and B. Mazur},
{\em Arithmetic Moduli of Elliptic Curves},
Ann.\ of Math.\ Stud.\ {\bf 108},
Princeton University Press, Princeton, 1985.

\bibitem{kenkumomose} {\sc M.\ A.\ Kenku and F.\ Momose},
Torsion points on elliptic curves defined over quadratic fields,
{\em Nagoya Math.\ J.} {\bf 109} (1988), 125--149.

\bibitem{kubert} {\sc D. Kubert},
Universal bounds on the torsion of elliptic curves,
{\em Proc.\ London Math.\ Soc.}, (3) {\bf 33} (1976), no.\ 2, 193--237.

\bibitem{kuhn}
{\sc R. M. Kuhn},
Curves of genus 2 with split Jacobian,
{\em Trans. Amer. Math. Soc.} {\bf 307} (1988), no.~1, 41--49.

\bibitem{leprevost1}
{\sc F. Lepr\'evost},
Famille  de courbes de genre $2$ munies d'une classe de diviseurs rationnels d'ordre $13$,
{\em C. R. Acad. Sci. Paris S\'er. I Math.} {\bf 313} (1991), no.~7, 451--454.

\bibitem{leprevost2}
{\sc F. Lepr\'evost},
Familles de courbes de genre $2$ munies d'une classe de diviseurs rationnels d'ordre
$15$, $17$, $19$ ou $21$,
{\em C. R. Acad. Sci. Paris S\'er. I Math.} {\bf 313} (1991), no.~11, 771--774.

\bibitem{leprevost3}
{\sc F. Lepr\'evost},
Torsion sur des familles de courbes de genre $g$,
{\em Manuscripta Math.} {\bf 75} (1992), no.~3, 303--326.

\bibitem{leprevost4}
{\sc F. Lepr\'evost},
Famille de courbes hyperelliptiques de genre $g$ munies d'une classe de diviseurs
rationnels d'ordre $2g^2 + 4g + 1$,
in:
{\em S\'eminaire de Th\'eorie des Nombres de Paris, 1991--1992},
Progress in Math. {\bf 116},
Birkh\"auser, Boston, (1993), 107--119.

\bibitem{leprevost5}
{\sc F. Lepr\'evost},
Points rationnels de torsion de jacobiennes de certaines courbes de genre $2$,
{\em C. R. Acad. Sci. Paris S\'er. I Math.} {\bf 316} (1993), no.~8, 819--821.

\bibitem{leprevost6}
{\sc F. Lepr\'evost},
Jacobiennes de certaines courbes de genre $2$ : torsion et simplicit\'e,
{\em J. Th\'eor. Nombres Bordeaux} {\bf 7} (1995), no.~1, 283--306.

\bibitem{leprevost7}
{\sc F. Lepr\'evost},
Sur une conjecture sur les points de torsion rationnels des jacobiennes de courbes,
{\em J. Reine Angew. Math.} {\bf 473} (1996), 59--68.

\bibitem{leprevost8}
{\sc F. Lepr\'evost},
Sur certains sous-groupes de torsion de jacobiennes de courbes hyperelliptiques
de genre $g\ge 1$, {\em Manuscripta Math.} {\bf 92} (1997), no.~1, 47--63.

\bibitem{leprevost9}
{\sc F. Lepr\'evost},
A provisional report on the rational torsion groups of jacobians 
of curves of genus two, preprint (1996).

\bibitem{mazur} {\sc B. Mazur}, Modular curves and the Eisenstein ideal,
{\em Inst. Hautes \'{E}tudes Sci. Publ. Math.} {\bf 47} (1977), 33--186 (1978).

\bibitem{mestreconducteurs} {\sc J.-F. Mestre},
Formules explicites et minorations de conducteurs de vari\'{e}t\'{e}s
alg\'{e}briques, {\em Compositio Math.} {\bf 58} (1986), no.~2, 209--232.

\bibitem{milne}
{\sc J. Milne},
Abelian varieties,
in: 
G. Cornell and J. H. Silverman (eds.),
{\em Arithmetic Geometry},
Springer-Verlag, New York, (1986), 103--150.


\bibitem{ogawa}
{\sc H. Ogawa},
Curves of genus~$2$ with a rational torsion divisor of order $23$,
{\em Proc. Japan Acad. Ser. A Math. Sci.} {\bf 70} (1994), no.~9, 295--298.

\bibitem{oortueno}
{\sc F.\ Oort and K. Ueno},
Principally polarized abelian varieties of dimension two or three are
Jacobian varieties,
{\em J. Fac. Sci. Univ. Tokyo Sect. IA Math.} {\bf 20} (1973), 377--381.

\bibitem{schaefer}
{\sc E.\ F.\ Schaefer},
$2$-descent on the Jacobians of hyperelliptic curves,
{\em J. Number Theory} {\bf 51} (1995), no.~2, 219--232.

\bibitem{silverman}
{\sc J.\ H.\ Silverman},
{\em The Arithmetic of Elliptic Curves},
Grad.\ Texts in Math.\ {\bf 106},
Springer-Verlag, New York, 1986.

\end{thebibliography}
\end{document}